\documentstyle[11pt,amsmath,amsfonts,amssymb]{report}

\newcommand{\Z}{\mbox{\msbm Z}}

\newfont{\msbm}{msbm9 at 11pt}
\newfont{\msbms}{msbm5 at 8pt}
\newcommand{\C}{\mbox{\msbm C}}

\newfont{\smsbm}{msbm10 at 9pt}
\newfont{\gmsbm}{msbm10 at 18pt}
\newfont{\Gmsbm}{msbm10 at 24pt}

\newcommand{\D}{\mbox{\msbm D}}

\newcommand{\N}{\mbox{\msbm N}}

\newcommand{\h}{\mbox{$\cal H$}}
\newcommand{\x}{\mbox{$\cal X$}}
\newcommand{\lx}{\mbox{$\cal L(X)$}}
\newcommand{\lh}{\mbox{$\cal L(H)$}}

\newcommand{\sst}{$\sigma (T)$}
\newcommand{\rrt}{$\rho (T)$}

\newcommand{\saa}{$\sigma_{ap}(T)$}

\newcommand{\lsx}{\sigma_{_T}}
\newcommand{\lrx}{\rho_{_T}}
\newcommand{\xtf}{{\cal X}_{_T}(F)}
\newcommand{\xxt}{{\cal X}_{_T}}
\begin{document}

\newtheorem{thm}{Theorem}[section]
\newtheorem{dft}[thm]{Definition}
\newtheorem{prop}[thm]{Proposition}
\newtheorem{lem}[thm]{Lemma}
\newtheorem{rem}[thm]{Remark}
\newtheorem{prf}{Proof}
\newtheorem{cor}[thm]{Corollary}
\thispagestyle{empty}
\vspace{2.00cm}
\begin{center}
{\LARGE\bf BOUNDED POINT EVALUATIONS AND }\\
\vspace{0.50cm}

{\LARGE\bf LOCAL SPECTRAL THEORY}\\
\vspace{3.00cm} 
{\Large \bf Abdellatif BOURHIM}\\
\vspace*{3.00cm}
{ The Abdus Salam International Centre for Theoretical Physics}\\
{Strada Costiera 11, Miramare P.O. Box 586, Trieste, Italy.}\\
\vspace{3.00cm} 
{\bf Abstract}
\end{center}
~~~~We study in this paper the concept of bounded point evaluations for cyclic operators. We give a negative answer to a 
question of L.R. Williams {\it Dynamic Systems and Apllications} 3(1994) 103-112. Furthermore, we generalize some results of Williams 
and give a simple proof of theorem 2.5 of \cite{33} that non normal hyponormal weighted shifts have fat local spectra.
\begin{center}
\vspace{1.00cm}
{Dissertation presented to ICTP Mathematics section}\\
{in Candidacy for the Diploma Degree}\\
\vspace{2.50cm}
August 2000\\
\vspace{2.50cm}
{Supervisor: Professor C.E. Chidume} \\
chidume@ictp.trieste.it
\end{center}
\newpage
\thispagestyle{empty}
\vspace*{8.00cm}

{\Huge {\bf To my loving father and mother}}\\
\\
\\
\\
\hspace*{3.00cm}{\Huge {\bf To my brothers and sisters.}}

\tableofcontents

\chapter*{Introduction}
\addcontentsline{toc}{chapter}{Introduction}
The main purpose of our dissertation is to study links between local spectral theory and the concept of bounded point evaluations for 
arbitrary cyclic operators on a Hilbert spaces.\\
\\
In chapter 1, we first fix some terminology and recall some basic notions concerning the spectral theory of a bounded operator
on complex Banach spaces and vector valued analytic functions. The concepts of spectrum and
resolvent set are introduced, the various subdivisions of the spectrum are studied and the analyticity of resolvent map 
on the resolvent set is proved.	The last part of this chapter is devoted to weighted shift operators which are widely
used for solving many problems in operator theory. Namely, we describe the spectrum and its parts of weighted shift operators.\\
\\    
In chapter 2 two basic concepts of local spectral theory are discussed. In section one are given basic properties of the operators with the single
valued extension property which have been considered by N. Dunford in 1952 those bounded operators $T$ on Banach
space ${\cal X}$ such that for every open set $U\subset \C$, the equation $(T-\lambda I)F(\lambda)=0$ admits the zero function
$F\equiv 0$ as a unique analytic solution. However, the generalization of the spectrum and resolvent set are introduced 
and studied. Finally in section two are studied the elementary properties of operators with Dunford's Condition C, which
will be used frequently throughout this dissertation.\\
\\
Chapter 3 is devoted to two classes of the most interesting operators which are defined around the notion of normal operators,
namely the subnormal and hyponormal operators. The concept of subnormal and hyponormal operators was introduced by Paul R. Halmos in 1950 
in \cite{222}, and many basic properties of subnormal operators were proved by Bram in \cite{5}. In this chapter, we present the basic tools for their understanding.
Are shown that the minimal normal extension of a subnormal operator $S$ on a Hilbert space $\h$ is unique up to an 
invertible isometry and its spectrum is contained in the spectrum of $S$ and the hyponormal operators have the single 
valued extension property and the Dunford's Condition C. Lastly, are given characterizations of normal, subnormal 
and hyponormal weighted shift operators.\\
\\
Chapter 4 is inspired by the papers of L.R. Williams \cite{32} and T.T. Trent \cite{300}. It contains a study of the bounded point evaluations
of cyclic operators and our results. We give a negative answer to the question A posed by Williams in \cite{32}, we compare
the concept of analytic bounded point evaluation of a weighted shift operators given by A.L. Shields in \cite{25} with the one
defined by L.R. Williams, we prove that if $T$ is a cyclic bounded operator on a Hilbert space $\h$ with Dunford's Condition C
and without point spectrum, then the local spectra
of $x$ with respect to $T$ is equal to the spectrum of $T$ for each $x$ in a dense subset of $\h$, finally, we give a simple
proof of the Theorem 2.5 of \cite{33} using the bounded point evaluation of weighted shifts and the fact that a non zero 
analytic function has isolate zeroes.    

\pagestyle{headings}
\chapter{Spectral Theory for Weighted Shift Operators}

\section{Preliminaries}
\hspace*{1.00cm}In this section we assemble some background material from spectral theory of operators 
which will be needed in the sequel.\\
\hspace*{1.00cm}In what follows, $\cal X^*$ will denote the dual space of a complex Banach space $\x$ and 
$\lx$ will denote the algebra of all linear bounded operators on $\x$. 
The spectrum \sst ~of an operator  $T\in \lx$ is defined as follows 
$$\sigma (T)=\big{\{}\lambda\in\C~:~ T-\lambda I\mbox{ is not invertible in }{\cal L(H)}\big{\}}.$$
It is a nonempty compact subset of $\C$ contained in the ball $\{\lambda \in \C
~~/~~|\lambda|\leq \|T\| \}$ (see \cite{13}). 
The spectral radius $r(T)$ of an operator $T\in \lx$ is defined by
$$r(T)=\sup\{|\lambda |~:~\lambda \in \sigma(T) \};$$
and satisfies 
$$r(T)=\lim_{n\to\infty}\|T^n\|^{\frac{1}{n}}\leq \|T\|$$
in the sense that the indicated limit always exist and has the indicated value (see \cite{15}).\\
If an operator $T\in \lx$ is invertible then 
$$\sigma (T^{-1})=\big{\{}\frac{1}{\lambda}~:~\lambda\in\sigma(T)\big{\}}$$ 
and so
$$\frac{1}{r(T^{-1})}=\inf\big{\{}|\lambda |~:~\lambda\in\sigma(T)\big{\}}$$
hence, in this case the spectrum  $\sigma(T)$ of $T$ is lies in the annulus $\big{\{}\lambda\in\C~:~\frac{1}{r(T^{-1})}\leq |\lambda |\leq
r(T)\big{\}}$.\\
It is often useful to ask of a point in the spectrum of an operator $T\in$\lx ~how it got there. The question
reduces therefore to this: why is a non-invertible operator not invertible? There are several possible ways of answering
the question; they have led to several classifications of spectra.      
Perhaps the simplest approach to the subject is to recall that an operator $T\in$\lx~
is bounded from below (i.e: $\inf\limits_{\|x\|=1}\|Tx\|>0$) and has a dense range if and only if it is invertible in \lx. Consequently,
$$\sigma (T)=\sigma_{ap}(T)\cup\Gamma (T)$$
where $\sigma_{ap}(T)$ is the approximate point spectrum of $T$, the set of complex numbers $\lambda$ such that
$T-\lambda I$ is not bounded from below, and $\Gamma (T)$ is the compression spectrum of $T$, the set 
of complex numbers $\lambda$ such that the range of $T-\lambda I$ is not dense in \x.\\
\\
\hspace*{1.00cm}The following results summarizes the main properties of the approximate point spectrum.
\begin{dft}The lower bound of an operator $T\in\lx$ is the value given by
$$m(T)=\inf\limits_{\|x\|=1}\|Tx\|.$$
\end{dft}
\begin{prop}Let $T$ be an operator in \lx. If $\lambda\in\C$ such that $|\lambda |<m(T)$ then $T-\lambda I$ is bounded below.
\end{prop}
{\bf Proof. } It suffice to observe that
$$m(T)-|\lambda |\leq \|Tx\|-|\lambda |\leq \|(T-\lambda I)x\|~~\mbox{for every }x\in{\cal X},~\|x\|=1.$$
$\hfill\blacksquare$
\begin{rem} Let $T$ and $S$ be two operators in \lx. Then:\\
{\rm(i)}~~~A complex number $\lambda$ is in $\sigma_{ap}(T)$ if and only if $m(T-\lambda I)=0.$\\
{\rm(ii)}~~
$$ m(T)\|x\|\leq\|Tx\|~~\mbox{for every }x\in{\cal X}.$$
{\rm(iii)}~If $T$ is invertible then $m(T)=\frac{1}{\|T^{-1}\|}.$\\
{\rm(iv)}$$m(T)m(S)\leq m(TS)\leq \|T\|m(S).$$
{\rm(v)}~$m(T)>0$ if and only if $T$ is one to one and has closed range.
\end{rem}
\begin{prop}Let $T\in\lx$. Then the sequence $\big{(}\big{(}m(T^n)\big{)}^{\frac{1}{n}}\big{)}_n$ is convergent and has a limit
equal to
$\sup\limits_{n\geq 1}\big{(}m(T^n)\big{)}^{\frac{1}{n}}.$ This limit will be denoted by $r_1(T).$
\end{prop}
{\bf Proof. }It is clear that 
$$\sup\limits_{n\geq 1}\big{(}m(T^n)\big{)}^{\frac{1}{n}}\geq\limsup\limits_{n\to \infty}\big{(}m(T^n)\big{)}^{\frac{1}{n}}.\leqno{(*)}$$
Now, fix $k$. Then for every $n$ there exists $p=p(n)$ and $r=r(n)$ such that $0\leq r<k$ and $n=kp+r.$ 
From Remark 1.1.3 it follows that 
$$m(T^n)\geq m(T^k)^pm(T)^r.$$ Hence $$\liminf\limits_n\big{(}m(T^n)\big{)}^{\frac{1}{n}}\geq
\big{(}m(T^k)\big{)}^{\frac{1}{k}}.$$
Since $k$ is arbitrary then
$$\liminf\limits_n\big{(}m(T^n)\big{)}^{\frac{1}{n}}\geq \sup\limits_n\big{(}m(T^n)\big{)}^{\frac{1}{n}}\leqno{(**)}$$
and the result follows from $(*)$ and $(**)$.$\hfill\blacksquare$
\begin{rem}The same argument of this proof can be used to prove that  
$$r(T)=\lim_{n\to\infty}\|T^n\|^{\frac{1}{n}}$$
for every operator $T$ in \lx.
\end{rem}
\begin{prop}For every operator $T$ in \lx,
$$\sigma_{ap}(T)\subset\big{\{}\lambda\in\C~:~r_1(T)\leq |\lambda |\leq r(T)\big{\}}.$$
\end{prop}
{\bf Proof. }Since $\sigma_{ap}(T)\subset\sigma (T)\subset\big{\{}\lambda\in\C~:~|\lambda |\leq r(T)\big{\}}$
then it suffice to prove that $r_1(T)\leq |\lambda |$ for every $\lambda\in\sigma_{ap}(T).$ So, let $\lambda\in\C$ such that
$|\lambda |<r_1(T).$ Then there $|\lambda |^n<m(T^n)$ for some integer $n$. By Proposition 1.1.2 we see that 
$$T^n-\lambda^nI=(T^{n-1}+\lambda T^{n-2}+...+\lambda^{n-1}I)(T-\lambda I)$$
is bounded below. Hence, it follows from the inequality (iv) in Remark 1.1.3 that $(T-\lambda I)$ is
bounded below; that is $\lambda\not\in\sigma_{ap}(T)$.$\hfill\blacksquare$
\begin{thm}Let $T\in$\lx.\\
{\rm(i)}~~~A complex number $\lambda$ belongs to $\sigma_{ap}(T)$ if and only if there exists a sequence $(x_n)$ of unit
vectors
in \x such that $\lim\limits_{n\to \infty}\|(T-\lambda I)x_n\|=0$. In particular,the approximate point spectrum of $T$ contains 
the set of eigenvalues of $T$, which is called the point spectrum of $T$ and denoted by $\sigma_p(T)$.\\
{\rm(ii)}~~The approximate point spectrum of $T$, \saa, is a closed non-empty subset of the spectrum of $T$.\\
{\rm(iii)}~The boundary of the spectrum of $T$ is contained in \saa.\\  
\end{thm} 
{\bf Proof. }The first property holds immediately from the definition of the approximate point spectrum of $T$.
Now we will show first that $\sigma_{ap}(T)\subset\sigma (T)$. Suppose $\lambda\not\in\sigma (T)$, then $(T-\lambda I)$ has bounded
inverse and so
$$M\leq \|(T-\lambda I)x\|~~~~{\mbox{for every }}x\in\x ~~~\|x\|=1$$
where $M=\frac{1}{\|(T-\lambda I)^{-1}\|}$, so $\lambda\not\in\sigma_{ap}(T)$. Hence $\sigma_{ap}(T)\subset\sigma (T)$.\\  
Let $\lambda\in \C\backslash\sigma_{ap}(T)$ then there exists $M>0$ such that 
$$M\leq \|(T-\lambda I)x\|~~~~{\mbox{for every }}x\in\x~~~\|x\|=1$$
Consequently, for every $\mu\in \C$ such that $|\lambda -\mu |<\frac{M}{2}$ we have
$$\frac{M}{2}\leq\|(T-\lambda I)x\|-|\lambda -\mu |\leq \|(T-\mu I)x\|~~~~{\mbox{for
every
}}x\in\x~~~\|x\|=1,$$
so $$\{\mu\in\C~~/~~|\lambda -\mu|<\frac{M}{2}\}\subset\C\backslash\sigma_{ap}(T).$$
Hence $\sigma_{ap}(T)$ is closed.\\
Next, let $\lambda$ be a point on the boundary of $\sigma (T)$ and let $\epsilon >0.$
Since $\sigma (T)$ and $\C\backslash\sigma (T)$ have the same boundary then there is $\mu\not\in \sigma (T)$ such that
$|\lambda -\mu |<\frac{\epsilon}{2}.$ We have 
$$\frac{2}{\epsilon}\leq \frac{1}{d(\mu )}\leq \|(T-\mu I)^{-1}\|~~\mbox{(see next section Theorem 1.2.1)}$$
where $d(\mu )$ denotes the distance from $\mu$ to the spectrum of $T$. Therefore there is $x\in$\x ~such that $\|x\|=1$ and
$$\frac{1}{\epsilon}\leq \|(T-\mu I)^{-1}x\|.$$
Let $y=\frac{1}{\|(T-\mu I)^{-1}x\|}(T-\mu I)^{-1}x$ then $\|y\|=1$ and 
$$\|(T-\lambda I)y-(T-\mu I)y\|<\frac{\epsilon}{2}.$$
It follows that
\begin{eqnarray*}
\|(T-\lambda I)y\|&\leq &\|(T-\lambda I)y-(T-\mu I)y\|+\|(T-\mu I)y\|\\
&<&\frac{\epsilon}{2}+\frac{\|x\|}{\|(T-\mu I)^{-1}\|}\\
&<&\frac{3\epsilon}{2}
\end{eqnarray*}
Since $\epsilon >0$ is arbitrary it follows that $\lambda\in\sigma_{ap}(T)$. Since a nonempty compact subset of \C~ has a nonempty
boundary then $\sigma_{ap}(T)$ is nonempty. The proof is complete.$\hfill\blacksquare$
\begin{thm}Let $T^*\in {\cal L(X^*)}$ be the adjoint operator of an bounded operator $T\in$\lx.~Then the following are equivalent:\\
{\rm(i)}~~~$T$ is bounded below.\\
{\rm(ii)}~~$T^*$ is surjective.
\end{thm}
For the proof, see \cite{4} Theorem 57.18 page 245.$\hfill\blacksquare$

\section{Analyticity of resolvents}
\hspace*{1.00cm}We first recall some basic notions and results from analyticity of functions whose domains are open sets in the complex plane and
whose values are vectors of the Banach space \x,~the monographs \cite{17}, \cite{18} and \cite{31} contain further informations.\\
\hspace*{1.00cm}Let $\Omega$ be a nonempty open subset of \C. A function $\rho$ from $\Omega$ to \x ~is said to be a {\bf analytic} in $\Omega$ if for every
$\phi\in\x^*$, the complex value function $\phi\circ\rho$ is analytic in the usual sense. In the case, when \x ~is a Hilbert space \h~ this
definition is equivalent to that for every $x\in$\h,~the function
\begin{eqnarray*}
&&\Omega\longrightarrow\C\\
&&\lambda\longmapsto\langle \rho (\lambda )~,~x\rangle
\end{eqnarray*}
is analytic in the usual sense. If \x ~ is the algebra \lh~of all bounded operators on a Hilbert space \h,~the analyticity of $\rho$ can be interpreted 
as following: For every $x,~y\in$\h ~the function 
\begin{eqnarray*}
&&\Omega\longrightarrow\C\\
&&\lambda\longmapsto\langle \rho (\lambda )x~,~y\rangle
\end{eqnarray*}
is analytic in the usual sense.\\
Using the fact that $0$ is the unique element $x\in$\x ~satiafaiying  $\phi (x)=0$ for every $\phi\in$\x $^*$, 
we  note that almost all the results concerning the analytic complex functions can be extended to vector
valued functions. Such an useful result is Liouville's Theorem which say that every bounded entire vector valued function 
is a constant.\\
\hspace*{1.00cm}Let $T\in$\lx ~ be a bounded operator on \x .~The complement set of the spectrum \sst ~of $T$, denoted by \rrt ,~is called the
resolvent set of $T$. It is a nonempty open subset of \C.~The mapping
\begin{eqnarray*}
&R_{_T}:&\rho (T)\longrightarrow{\cal L(H)}\\
&&\lambda\longmapsto (T-\lambda I)^{-1}
\end{eqnarray*}
is called the resolvent map of $T$. 
\begin{thm} Let $T\in$\lx.~The resolvent map of $T$ is analytic on \rrt~vanishing at infinity and has the following property
$$\frac{1}{d(\lambda )}\leq \|R_{_T}(\lambda )\|~~for ~every~\lambda\in\rho (T) $$
where $d(\lambda )$ is the distance from $\lambda$ to the spectrum \sst~of $T$. Therefore $\|R_{_T}(\lambda )\|\to\infty$ as $d(\lambda )\to 0.$
\end{thm}
\hspace*{1.00cm}In the proof of this Theorem we shall require the following Lemma:
\begin{lem}Let $T\in\lx$ be a bounded operator on \x~.If $\|T\|<1$ then $I+T$ is invertible and
$$\|(T-I)^{-1}+I+T\|\leq \frac{\|T\|^2}{1-\|T\|}.$$
In particular, if $\lambda$ is a fixed point in \rrt~then for every $\mu\in\C$ with $|\mu |<\|(T-\lambda I)^{-1}\|^{-1}$ we have $\lambda
+\mu\in\rho (T).$
\end{lem}
{\bf Proof. }Since $\|T\|<1$ and $\|T^k\|\leq \|T\|^k$ for every $k\geq 0$ then this series $\sum\limits_{n=0}^{+\infty}T^n$ converges
in \lx~and,
$$(T-I)\bigg{(}\sum\limits_{n=0}^{+\infty}T^n\bigg{)}=\bigg{(}\sum\limits_{n=0}^{+\infty}T^n\bigg{)}(T- I)=I.$$
Hence, $T- I$ is invertible and
$$(T-I)^{-1}=\sum\limits_{n=0}^{+\infty}T^n.$$ 
The rest of this proof follows easily.\\
{\bf Proof of Theorem 1.2.1. }Let $\lambda$ be a fixed point in \rrt~. Then there is $0<r$ such that $\mu\in$\rrt~for
every $\mu\in\C$ with $|\lambda -\mu |<r$. Since
$$\frac{(T-\lambda I)^{-1}-(T-\mu I)^{-1}}{\lambda -\mu}=(T-\lambda I)^{-1}(T-\mu I)^{-1},$$
then
$$\lim\limits_{\mu\to\lambda}\frac{(T-\mu I)^{-1}-(T-\lambda I)^{-1}}{\mu -\lambda}=(T-\lambda I)^{-2}.$$
Therefore, for every $\phi\in$\lx$^*$ 
$$\lim\limits_{\mu\to\lambda}\frac{\phi\circ R_{_T}(\mu )-\phi\circ R_{_T}(\lambda )}{\mu -\lambda}=\phi\circ(T-\lambda
I)^{-2}.$$ 
Hence, $\phi\circ R_{_T}$ is differentiable on \rrt.~So, $R_{_T}$ is analytic on \rrt.
On the other hand, we have
$$\lambda R_{_T}(\lambda )=(\frac{1}{\lambda}T-I)^{-1},$$
this shows that the resolvent $R_{_T}$ of $T$ vanishing at infinity.\\
Since \sst~ is a nonempty compact set then there is $\mu_0\in$\sst~such that $d(\lambda )=|\mu_0-\lambda |$. So, if
$d(\lambda )=|\mu_0-\lambda |<\|R_{_T}(\lambda )\|^{-1}=\|(T-\lambda I)^{-1}\|^{-1}$ then it follows from Lemma 1.2.2 that
$\mu_o=(\mu_0-\lambda )+\lambda\in $\rrt~. Contradiction. This completes the proof of the Theorem.$\hfill\blacksquare$


\section{Weighted Shift Operators}
\hspace*{1.00cm}A {\bf weighted shift operator $T$ on complex Hilbert space }\h ~is an operator that maps each vector in some
orthonormal basis $(e_n)_n$
into a scalar multiple of the next vector
$$Te_n=\omega_n e_{n+1}$$
for all $n$. The operator $T$ is called a unilateral weighted shift when the index $n$ runs the nonnegative integers and it is
called a bilateral weighted shift when the $n$ runs over all integers.\\
In what follows, $T$ will always denote a weighted shift operator with a weight sequence $(\omega_n)_n$. We shall
sometimes omit the adjective "weighted" and refer to $T$ simply as a shift.\\
\hspace*{1.00cm}An operator $A$ on a $n$-dimensional Hilbert space is called a {\bf finite-dimensional weighted shift}    
if there are numbers $\{\alpha_1,...,\alpha_{n-1}\}$ and an orthonormal basis $v_1,...,v_n$ such that
\begin{eqnarray*}
Av_k&=&\alpha_kv_{k+1}~~~(k<n)\\
Av_n&=&0.
\end{eqnarray*}
Such operator is nilpotent (i.e: there is a positive integer $s$ such that $A^s=0$).\\
Note that a weighted shift operator is injective if and only if none of the weights is zero. If finitely
many weights are zero then $T$ is a direct sum of a finite -dimensional weighted shifts and an 
infinite-dimensional injective weighted shift. If infinitely many weights are zero then $T$ is the direct
sum of an infinite family of finite-dimensional weighted shifts. This situation can be used to give an
example of an operator with spectral radius $1$, which is the direct sum of countable family of finite
dimensional nilpotent operators each of which has spectrum $\{0\}$ (see \cite{15}).\\
From now we shall assume that none of the weights is zero and let $\beta $ be the following sequence given
by:  

\begin{displaymath}
\beta_n=\left\{
\begin{array}{ll}
\omega_0...\omega_{n-1}&\textrm{\mbox{if }$n>0$}\\
\\
1&\textrm{\mbox{if }$n=0$}\\
\\
\frac{1}{\omega_{n}...\omega_{-1}}&\textrm{\mbox{if }$n<0$}\\
\end{array}
\right.
\end{displaymath}
\hspace*{1.00cm}An operator $U$ in \lh ~is called {\bf unitary} if $UU^*=U^*U=I$ which is equivalent that $U$ 
is invertible and $U^{-1}=U^*$. Two bounded operators $A$ and $B$ on \h ~are said to be
{\bf unitarily equivalent} if there is unitary operator $U$ in \lh ~such that $AU=UB.$ Therefore, two unitarily
equivalent bounded operators on \h ~have the same spectrum, the same point spectrum, the
same approximate point spectrum and the same compression spectrum.
\subsection{Elementary properties}
\begin{prop}
The shift $T$ is bounded if and only if the weight sequence is bounded. In this case,
\begin{eqnarray*}
\|T^n\|&=&\sup\limits_k|w_kw_{k+1}...w_{n+k-1}|\\
&=&\sup\limits_k\big{|}\frac{\beta_{n+k}}{\beta_k}\big{|}.
\end{eqnarray*}
\end{prop}
{\bf Proof. }
First suppose that the shift $T$ is bounded then for every integer $n$ we have
$$\|Te_n\|=|\omega_n|\leq\|T\|$$
so the weight $(\omega_n)_n$ is obviously bounded. Conversely, suppose there is $M>0$ such that
$|\omega_n|\leq M \mbox{ for every }n$ then for every $x=\sum\limits_na_ne_n\in\h$ we have
$\|Tx\|^2=\sum\limits_n|a_n\omega_n|^2$ therefore
$$\|Tx\|\leq M\big{(}\sum\limits_n|a_n|^2\big{)}^{\frac{1}{2}}=M\|x\|\mbox{  for every }x\in {\cal H}.$$
Thus $T$ is bounded. The equalities follow from the following relation 
$$T^ne_k=\omega_k\omega_{k+1}...\omega_{n+k-1}e_{n+k}.$$
$\hfill\blacksquare$\\
\hspace*{1.00cm}In the sequel we suppose that the weights $(\omega_n)_n$ is bounded just to have that $T$ is in \lh.
\begin{prop}If $T$ is bilateral weighted shift then 
$$T^*e_n=\overline{\omega_{n-1}}e_{n-1} \mbox{ for every }n.$$
If $T$ is unilateral weighted shift then
\begin{eqnarray*}
&&T^*e_n=\overline{\omega_{n-1}}e_{n-1} \mbox{ for every }n\geq 1\\
&&T^*e_0=0
\end{eqnarray*}
\end{prop}
{\bf Proof. } For every integers $n$ and $k$ we have
\begin{eqnarray*}
\langle T^*e_n~,~e_k\rangle &=&\langle e_n~,~Te_k\rangle \\
&=&\langle e_n~,~\omega_ke_{k+1}\rangle \\
&=&\overline{\omega_k}\langle e_n~,~e_{k+1}\rangle \\
\end{eqnarray*}
and the result follows.$\hfill\blacksquare$
\begin{rem} It follows from Proposition 1.3.2 that the unilateral shift $T$ is never invertible because 
$T^*$ is not, but the bilateral shift $T$ can be invertible. It is invertible if and only if the weight
$(\frac{1}{\omega_n})_n$ is bounded if and only if $\inf \frac{\beta_{n+1}}{\beta_n}>0$. In this case, for
$n=0,1,2,...$
\begin{eqnarray*}
\|T^{-n}\|&=&\sup\limits_k\frac{1}{|w_kw_{k+1}...w_{n+k-1}|}\\
&=&\sup\limits_k\big{|}\frac{\beta_k}{\beta_{n+k}}\big{|}\\
&=&\big{[}\inf\limits_k\frac{\beta_{n+k}}{\beta_k}\big{]}^{-1}.
\end{eqnarray*}
For this, it suffices to observe that $T^{-1}v_k=\frac{1}{\omega_{-k}}v_{k+1}$ where $v_k=e_{-k}$. Thus
$T^{-1}$ is a bilateral weighted shift and apply Proposition 1.3.1. 
\end{rem}
\begin{prop}If $(\lambda_n)$ are complex numbers of modulus 1, then $T$ is unitarily equivalent to the weighted shift
operator with weight sequence $(\overline{{\lambda}_{n+1}}\lambda_n\omega_n)_n$. 
\end{prop}
{\bf Proof. }Let $U$ be the unitarily operator defined by $Ue_n=\lambda_ne_n$. Then the operator $U^*TU$ is a
weighted shift with the weight sequence indicated above.$\hfill\blacksquare$
\begin{cor}
$T$ is unitarily equivalent to the weighted shift operator with weight sequence $(|\omega_n|)_n.$
\end{cor}
{\bf Proof. }Choose $\lambda_0=1$; the choice of the remaining $\lambda_n$ is then forced in both the unilateral and the
bilateral cases.$\hfill\blacksquare$
\begin{cor}
If $|c|=1$, then $T$ and $cT$ are unitarily equivalent.
\end{cor}
{\bf Proof. }Take $\lambda_n=\bar{c}^n$ for every $n$ and apply Proposition 1.3.5.$\hfill\blacksquare$
\begin{rem}In this section, our goal is to describe the spectrum, the point spectrum and the approximate
point spectrum of the weighted shift operators. We have already mentioned that two unitarily
equivalent bounded operators on \h ~have the same spectrum, the same point spectrum and the
same approximate point spectrum. Therefore by the Corollary 1.3.5, we may assume that the weights are
nonnegative and by the Corollary 1.3.6 we see that the spectrum, the point spectrum and the approximate
point spectrum of $T$ have circular symmetry about the origin. In particular, the entire circle
$\big{\{}\lambda\in\C~~/~~|\lambda |=r(T)\big{\}}$ is in the spectrum \sst ~of $T$. When $T$ is bilateral invertible shift then by
the same reason as before, the entire circle $\big{\{}\lambda\in\C~~/~~|\lambda |=\frac{1}{r(T^{-1})}\big{\}}$ is also in the
spectrum \sst ~of $T$.
\end{rem}
\subsection{Spectrum of Weighted Shifts}

\begin{thm}
If $T$ is unilateral shift then its spectrum is the disc
$$\{\lambda\in \C~~/~~|\lambda |\leq r(T)\}.$$
\end{thm}
{\bf Proof. } It is known that
$$\sigma (T)\subset\{\lambda\in \C~~/~~|\lambda |\leq r(T)\}.$$
Conversely, let $\lambda $ be in the resolvent set \rrt~ of $T$ and set $x=(T-\lambda I)^{-1}e_0=\sum\limits_n a_ne_n$,
then 
$$ a_{n}=-\frac{1}{\lambda^{n+1}}\beta_{n} {\mbox{~~for every }}n.$$
On the other hand, we have

\begin{eqnarray*}
\langle (T-\lambda I)^{-1}e_n~,~e_{n+k}\rangle &=&\frac{1}{\beta_n}\langle (T-\lambda
I)^{-1}T^ne_0~,~e_{n+k}\rangle\\
&=&\frac{1}{\beta_n}\langle T^n(T-\lambda I)^{-1}e_0~,~e_{n+k}\rangle\\ 
&=&\frac{1}{\beta_n}\langle T^nx~,~e_{n+k}\rangle\\
&=&\frac{\beta_{n+k}}{\beta_n\beta_k}a_k\\
&=&-\frac{\beta_{n+k}}{\beta_n}\frac{1}{\lambda^{k+1}}.
\end{eqnarray*}
Using Cauchy-Schwartz inequality we get,
$$\frac{\beta_{n+k}}{\beta_n}\frac{1}{|\lambda |^{k+1}}\leq \|(T-\lambda I)^{-1}\|.$$
By passing to the supremum on $n$ we get 
$$\|T^k\|\leq |\lambda |^{k+1}\|(T-\lambda I)^{-1}\|.$$
Taking $k$th roots and letting $k\to\infty$ gives $r(T)\leq |\lambda |$. Strict inequality must
hold since the entire circle $\{\lambda\in\C~~/~~|\lambda |=r(T)\}$ is in the spectrum of $T$ (see Remark 1.3.7). 
This proves the theorem.$\hfill\blacksquare$
\begin{thm}
{\rm (i)}~~~If $T$ is a invertible bilateral weighted shift then its spectrum is the annulus 
$$\big{\{}\lambda\in\C~~/~~\frac{1}{r(T^{-1})}\leq |\lambda|\leq r(T)\big{\}}.$$
{\rm (ii)}~~If $T$ is a bilateral weighted shift that is not invertible then its spectrum is the disc
$$\{\lambda\in \C~~/~~|\lambda |\leq r(T)\}.$$   
\end{thm}
{\bf Proof. }Let $\lambda\not= 0$ lie in the resolvent set
$\rho (T)$ of $T$ then for $x=\sum\limits_{n\in\Z}a_ne_n=((T-\lambda I)^{-1})e_0$ we have
$$
\leqno{(*)}~~~~~~~~~~~~~~~~~~~~~~~~~~~~~~~~~\left\{
\begin{array}{lll}
a_{-1}\omega_{-1}-\lambda a_0=1\\
\\
a_n\omega_n-\lambda a_{n+1}=0 ~~~\mbox{  for every } n\not= -1
\end{array}
\right.
$$
and so,
$$a_n=\frac{1}{\lambda^n}a_0\beta_n~~\mbox{and  }a_{-n-1}=\lambda^n\beta_{-n-1}\omega_{-1}a_{-1}~~~\mbox{for every }n>0.\leqno{(**)}$$
By the same calculation as in the proof of Theorem 1.3.8 we have
$$|a_k|\leq \frac{\beta_n\beta_k}{\beta_{n+k}}\|(T-\lambda I)^{-1}\|~~\mbox{for every }n\in\N ~\mbox{and }k\in\Z.\leqno{(***)}$$
By $(*)$ either $a_0$ or $a_1$ is non zero. Consider the following two cases.\\
{\bf Case 1:} $a_0\not= 0$. If we multiply $(***)$ by $|\lambda |^k$ and apply $(**)$ we obtain
$$|a_0|\frac{\beta_{n+k}}{\beta_n}\leq |\lambda|^k\|(T-\lambda I)^{-1}\|~~\mbox{for every }k\in\N.$$
Hence by the same argument of the proof of Theorem 1.3.8 we get $r(T)<|\lambda|.$\\
{\bf Case 2:} $a_1\not= 0$. In $(***)$ we take $k=-m \mbox{  for }m\geq 1$ and by applying the second equation in $(**)$ we get
$$|\lambda |^{m-1}\frac{\beta_{n-m}}{\beta_n}\leq\frac{1}{\omega_{-1}a_{-1}}\|(T-\lambda I)^{-1}\|.$$
In this case $T$ is invertible (just take in the last inequality $m=1$). By passing to the infimum on $n$ we get
$$|\lambda |^{m-1}\|T^{-m}\|\leq\frac{1}{\omega_{-1}a_{-1}}\|(T-\lambda I)^{-1}\|.$$
Taking $m$th roots and letting $m\to\infty$ we obtain $|\lambda |\leq \frac{1}{r(T^{-1})}$.
The equality is excluded by circular symmetry about the origin. So, $|\lambda |<\frac{1}{r(T^{-1})}$.\\
Since $a_{-1}\omega_{-1}-\lambda a_0=1$, then at least one of these two cases must hold for each given $\lambda$. 
On the other hand they cannot both occur since the conclusions exclude one another. If $T$ is invertible then the two cases together yield
the desired conclusion. If $T$ is not invertible then $\inf\limits_n\frac{\beta_{n+1}}{\beta_n}=0.$ So, Taking $k=-1$ in the inequality
$(***)$ and by passing to the infimum on $n$ the second case will be excluded; hence the first case yields the desired conclusion.$\hfill\blacksquare$ 
\begin{thm}If $T$ is unilateral weighted shift then
$$\sigma_{ap}(T)=\big{\{} \lambda\in\C~:~r_1(T)\leq |\lambda |\leq r(T)\big{\}}.$$
\end{thm}
{\bf Proof. }It is known that $\sigma_{ap}(T)\subset\big{\{}\lambda\in\C~:~r_1(T)\leq |\lambda |\leq r(T)\big{\}}$ (see Proposition 1.1.6), so if
$r_1(T)=r(T)$ then by nonemptiness and circular symmetry of the approximate point spectrum of $T$ we have
$\sigma_{ap}(T)=\big{\{}\lambda\in\C~:~|\lambda |=r(T)\big{\}}.$ Now suppose that $r_1(T)<r(T)$. Since $\sigma_{ap}(T)$ is
closed and has circular symmetry it suffice to prove that every positive real $c,~r_1(T)<c<r(T)$ is lies in $\sigma_{ap}(T).$
Let $c\in \big{(}r_1(T)~,~r(T)\big{)}$ and choose two reals numbers $a,~b$ such that $r_1(T)<a<c<b<r(T).$ Let $\epsilon >0,$
since $r(T)=\lim\limits_{n\to\infty}\big{[}\sup\limits_k\frac{\beta_{n+k}}{\beta_k}\big{]}^{\frac{1}{n}}$ and
$\frac{c}{b}<1$ then there exists $n$ and $k$ such that
$$\big{[}{\frac{c}{b}}\big{]}^n<\epsilon ~~\mbox{and  }\big{[}\frac{\beta_{n+k}}{\beta_k}\big{]}^{\frac{1}{n}}>b.$$
Also, we have $r_1(T)=\lim\limits_{m\to\infty}\big{[}\inf\limits_p\frac{\beta_{m+p}}{\beta_p}\big{]}^{\frac{1}{m}}$ 
and $\frac{a}{c}<1$ then there exist $p$ and $m>n+k$ such that
$$\big{[}{\frac{a}{c}}\big{]}^p<\epsilon ~~\mbox{and  }\big{[}\frac{\beta_{m+p}}{\beta_p}\big{]}^{\frac{1}{p}}<a.$$
Let $x=\sum\limits_sa_se_s\in \h $ where 
\begin{eqnarray*}
&&a_{k}=1,\\
&&a_s=\frac{\beta_{s}}{\beta_kc^{s-k}}~~\mbox{if } k+1\leq s\leq m+p+1,\\
&&a_s=0~~\mbox{otherwise.}
\end{eqnarray*}
And so,
$$a_{s-1}\omega_{s-1}-ca_s=0~~\mbox{for every  }k<s<m+p,\leqno{(*)}$$
$$\frac{1}{a_{n+k}}=\frac{\beta_kc^n}{\beta_{n+k}}<\big{[}{\frac{c}{b}}\big{]}^n<\epsilon ,\leqno{(**)}$$
and
$$\frac{a_{m+k}}{a_m}=\frac{\beta_{m+p}}{\beta_mc^p}<\big{[}{\frac{a}{c}}\big{]}^p<\epsilon .\leqno{(***)}$$
Therefore from $(*)$ it follows
$$Tx-cx=a_{m+p}\omega_{m+p}e_{m+p+1}-ca_ke_k$$
and hence
\begin{eqnarray*}
\|Tx-cx\|^2&=&{|a_{m+p}\omega_{m+p}|}^2+c^2\\
&\leq &\|T\|^2\big{(}|a_{m+p}|^2+1\big{)}
\end{eqnarray*}
because $c<r(T)\leq \|T\|=\sup\limits_r\omega_r.$
On the other hand, we have
$$\|x\|^2=\sum\limits_r|a_r|^2\geq |a_m|^2+|a_{n+k}|^2$$
then
\begin{eqnarray*}
\frac{\|Tx-cx\|^2}{\|x\|^2}&\leq &\|T\|^2\frac{|a_{m+p}|^2+1}{|a_m|^2+|a_{n+k}|^2}\\
&\leq &\|T\|^2\max\big{(}\frac{1}{|a_{n+k}|^2}~,~\frac{|a_{m+p}|^2}{|a_m|^2}\big{)}\\
&<&\epsilon^2\|T\|^2~~~\mbox{by }(**) ~~\mbox{and } (***)
\end{eqnarray*}
and so, $c$ is in $\sigma_{ap}(T)$.$\hfill\blacksquare$
\\
\hspace*{1.00cm}In dealing with bilateral shifts we shall use the following notations.
$$r_1^+=\lim\limits_{n\to +\infty}\big{[}\inf\limits_{j\geq 0}\frac{\beta_{n+j}}{\beta_j}
\big{]}^{\frac{1}{n}},~~~r^+=\lim\limits_{n\to +\infty}\big{[}\sup\limits_{j\geq 0}
\frac{\beta_{n+j}}{\beta_j}\big{]}^{\frac{1}{n}}$$
$$r_1^-=\lim\limits_{n\to +\infty}\big{[}\inf\limits_{j<0}\frac{\beta_{j}}{\beta_{-n+j}}
\big{]}^{\frac{1}{n}},~~~r^-=\lim\limits_{n\to +\infty}\big{[}\sup\limits_{j<0}
\frac{\beta_{j}}{\beta_{-n+j}}\big{]}^{\frac{1}{n}}$$
\\
\hspace*{1.00cm}The proof of the following theorem require a lot of technical computations, we omit it here.
\begin{thm}{\rm (RIDGE \cite{24})}. If $T$ is bilateral shift and if $r^-<r_1^+,$then
$$\sigma_{ap}(T)=\big{\{}\lambda\in\C~:~r_1^-\leq |\lambda|\leq
r^-\big{\}}\cup\big{\{}\lambda\in\C~:~r_1^+\leq |\lambda|\leq r^+\big{\}}.$$
Otherwise
$$\sigma_{ap}(T)=\sigma(T)=\big{\{}\lambda\in\C~:~\min (r_1^-,r_1^+)\leq |\lambda |\leq
\max (r^-,r^+)\big{\}}.$$
\end{thm}
\begin{thm}If $T$ is unilateral shift then:\\
{\rm (i)}~~~$\sigma_p(T)$ is empty.\\
{\rm (ii)}~~$\sigma_p(T^*)=\{0\}$ if $r_2(T)=0$ otherwise
$$\big{\{}\lambda\in\C~:~|\lambda |<r_2(T)\big{\}}\subset\sigma_p(T^*)\subset\big{\{}\lambda\in\C~:~|\lambda
|\leq r_2(T)\big{\}}$$
where $r_2(T)=\liminf\limits_{{n\to +\infty}}(\beta_n)^{\frac{1}{n}}.$ Furthermore, all eigenvalues of $T^*$ are simples.
\end{thm}
{\bf Proof. }(i) let $\lambda\in\C$ and let $x=\sum\limits_{n\in\N}a_ne_n\in$\h~ such that $Tx=\lambda x$. Then 
$$\lambda a_0=0~~\mbox{and }a_n\omega_n=\lambda a_{n+1} ~~~\mbox{for every }n\in\N.$$
Since $T$ is injective then $\lambda\not= 0$. Therefore, $a_n=0$ for every $n\in\N$. So, $x=0$. Thus $\sigma_p(T)$ is empty.\\
(ii) We have $T^*e_0=0$ then $0\in \sigma_p(T^*).$ Let $\lambda\in\sigma_p(T^*)$ and let $x=\sum\limits_{n\in\N}a_ne_n$ be
a corresponding eigenvector that is $T^*x=\lambda x.$ And so, $\lambda a_n=\omega_na_{n+1}$ for every $n\in\N$.
Therefore,
$$a_n=\frac{a_0\lambda^n}{\beta_n} ~~\mbox{for every } n\geq 1.$$
Hence $a_0\not= 0$ and $x=a_0\big{(}e_0+\sum\limits_{n\geq 1}\frac{\lambda^n}{\beta_n}e_n\big{)}.$
Therefore the eigenvalues of $T^*$ are simple and 
$$\|x\|^2=|a_0|^2\big{(}1+\sum\limits_{n\geq 1}\frac{|\lambda |^{2n}}{\beta_n^{2}}\big{)}.$$
By the Cauchy-Hadamard formula for the radius of convergence we get that $|\lambda |\leq r_2(T)$.
So,$$\sigma_p(T^*)\subset\big{\{}\lambda\in\C~:~|\lambda |\leq r_2(T)\big{\}}.$$
On the other hand, if $\lambda <r_2(T)$ then the series $\sum\limits_{n\geq 1}\frac{|\lambda |^{2n}}{\beta_n^{2}}$
is convergent and so, $x=e_0+\sum\limits_{n\geq 1}\frac{\lambda^n}{\beta_n}e_n$ is an eigenvector 
corresponding to $\lambda$. Thus $$\big{\{}\lambda\in\C~:~|\lambda |<r_2(T)\big{\}}\subset\sigma_p(T^*).$$
This proves the desired result.$\hfill\blacksquare$
\begin{rem}By circular symmetry, one of the two containments in (ii) must be equality.\\
If the unilateral weighted shift $T$ is not injective then $\sigma_p(T)=\{0\}$.
\end{rem}
\hspace*{1.00cm}If $T$ is bilateral shift, set:\\ 
$$r_2^+=\liminf\limits_{{n\to +\infty}}(\beta_n)^{\frac{1}{n}},~~~r_3^+=
\limsup\limits_{{n\to +\infty}}(\beta_n)^{\frac{1}{n}}$$
$$r_2^-=\liminf\limits_{{n\to +\infty}}\big{(}{\frac{1}{\beta_{-n}}}\big{)}^{\frac{1}{n}},~~~r_3^+=
\limsup\limits_{{n\to +\infty}}\big{(}{\frac{1}{\beta_{-n}}}\big{)}^{\frac{1}{n}}.$$
Clearly we have
$$r_1^-\leq r_2^-\leq r_3^-\leq r^-,~~~r_1^+\leq r_2^+\leq r_3^+\leq r^+$$   
   
\begin{thm}If $T$ is bilateral shift then:\\
{\rm (i)}~~~All eigenvalues of $T$ and $T^*$ are simple.\\
{\rm (ii)}~~$\big{\{}\lambda\in\C~:~r_3^+<|\lambda |<r_2^-\big{\}}\subset\sigma_p(T)\subset\big{\{}\lambda\in\C~:~r_3^+\leq |\lambda
|\leq r_2^-\big{\}}.$\\
{\rm (iii)}~$\big{\{}\lambda\in\C~:~r_3^-<|\lambda |<r_2^+\big{\}}\subset\sigma_p(T^*)\subset\big{\{}\lambda\in\C~:~r_3^-\leq |\lambda
|\leq r_2^+\big{\}}.$\\
{\rm (iv)}~~At least one of $\sigma_p(T),~\sigma_p(T^*)$ is empty. 
\end{thm}
{\bf Proof. }Let $\lambda\in\sigma_p(T)$ and let $x=\sum\limits_{n\in\Z}a_ne_n$ be an eigenvector corresponding to $\lambda.$ Then
$$a_{n-1}\omega_{n-1}=\lambda a_n~\mbox{for every }n\in\Z.$$
And so,
$$a_n=\frac{a_0\lambda^n}{\beta_n},~~a_{-n}=a_0\beta_{-n}\lambda^n ~~\mbox{for every } n\geq 1.$$
From this we wee that the eigenvalues are simple. Further, $\lambda$ is eigenvalue for $T$ if and only if the sequence $(a_n)_{n\in\Z}$
defined above is square summable. This leads to two power series, one in $\lambda$ and the other in $\frac{1}{\lambda}$, and the result follows
from the formula for the radius of convergence.\\
The case of $T^*$ is similar.\\  
Finally, assume that $\lambda$ and $\mu$ are eigenvalues of $T$ and $T^*$ respectively. We wish to show that this is impossible. By what
has gone before we must have $$r_3^+\leq |\lambda |\leq r_2^-~~~\mbox{and}~~~r_3^-\leq |\mu |\leq r_2^+.$$
Since $r_2^-\leq r_3^-$ and  $r_2^+\leq r_3^+,$ then $|\lambda |=|\mu |.$ Also, an examination of the series which must converge shows that
$$\sum\limits_{n\geq 1}\frac{|\mu |^{2n}}{\beta_n^2}<\infty ~~~\mbox{and}~~~\sum\limits_{n\geq 1}\frac{\beta_n^2}{|\lambda
|^{2n}}<\infty$$
which is impossible since $|\lambda |=|\mu |$.$\hfill\blacksquare$
\begin{rem} ~~\\
{\rm (i)}~~~If $r_2^-<r_3^+$ then $\sigma_p(T)$ is empty; if $r_2^+<r_3^-$ then $\sigma_p(T^*)$ is empty.\\
{\rm (ii)}~~By circular symmetry, one of the containments in {\rm (ii)}, and one of the containments in {\rm (iii)} must be equality.\\
{\rm (iii)}~Let $\omega_n=1$ for every $n\in$\Z. Then $T$ is unitary bilateral shift and its adjoint which is its inverse is given by
$T^*e_n=T^{-1}e_n=e_{n-1}$ for every $n$. 
On the other hand $T$ and $T^*$ are unitarily equivalent since $R^{-1}TR=T^*$ where $R$ is the unitary operator determined by
$Re_n=e_{-n}$ for every $n\in$\Z. So, $T$ and $T^*$ have the same spectrum, the same point spectrum and the same approximate point
spectrum. Hence from the property (iv) it follows that $\sigma_p(T)=\sigma_p(T^*)=\emptyset.$\\
Also, in this case we have $r(T)=r_1(T)=r_i^+=r_i^-=1,~~i=1,~2,~3.$ So,
$$\sigma (T)=\sigma_{ap}(T)=\sigma (T^{-1})=\sigma_{ap}(T^{-1})=\big{\{}\lambda\in\C~:~|\lambda |=1\big{\}}.$$ 
\end{rem}
\chapter{Analytic Extension and Spectral Theory}

\section{Operators with the single-valued Extension Property}
\hspace*{1.00cm}In this section, we introduce a generalization of the spectrum and the resolvent set due to Dunford (see \cite{2}).
Let \x~be a complex Banach space and \lx~be the algebra of all bounded operators on \x.~For an operator
$T\in$\lx~, the resolvent map $R_{_T}$ is a operator valued analytic function on the resolvent set
\rrt~of $T$. For $x\in\x$ the vector valued function $R_{_{T,x}}$ defined on \rrt~by
$$R_{_{T,x}}(\lambda )=R_{_T}(\lambda )x=(T-\lambda I)^{-1}x$$ is analytic function on \rrt~ and
satisfying the following equation
$$(T-\lambda I)R_{_{T,x}}(\lambda )=x~~\mbox{for every }\lambda\in\rho(T).$$
In many cases the function $R_{_{T,x}}$ can be extended to be analytic on an open set properly
containing \rrt.~We will call a vector valued analytic function $F$ an analytic extension of
$R_{_{T,x}}$ if the domain $D(F)$ of $F$ is containing \rrt~and
$$(T-\lambda I)F(\lambda )=x~~~\mbox{for every }\lambda\in D(F).$$
We now encounter the possibility that there may be many extensions of $R_{_{T,x}}$ and they may not
agree on their common domain. However if all extensions do agree on their common domain for each
$x\in\x$ we say that the operator $T$ has the {\bf single valued extension property} (s.v.e.p). Therefore,
for such operator $T$ and for every $x\in\x$, the function $R_{_{T,x}}$ has a maximal extension called the maximal single
valued extension of $R_{_{T,x}}$ and denoted by $\widetilde{x}(.)$. Thus
$\widetilde{x}(.)$ is an analytic vector valued function such that
$$(T-\lambda I)\widetilde{x}(\lambda )=x~~~\mbox{for every }\lambda\in D(\widetilde{x});$$
the open domain $D(\widetilde{x})$ of $\widetilde{x}$ will be called the {\bf local resolvent set} of $x$ and is denoted by
$\rho_{_T}(x)$ and its complement, denoted by $\sigma_{_T}(x)$, will be called the {\bf local spectrum} of $x$; it is 
a closed subset of the spectrum \sst~of $T$. We will see by using Liouville's Theorem that the local 
spectrum of $x$ is empty if and only if $x$ is the zero vector.\\
\hspace*{1.00cm}Before proceeding further it would be well to note that, an operator $T\in\lx$ has the single valued extension property (s.v.e.p) if and
only if  for every open set $U\subset\C$, the only analytic solution of the equation $(T-\lambda I)F(\lambda )=0~~\mbox{for }\lambda\in U$
is the zero function $F\equiv 0.$ Therefore, if $T\in\lx$ has (s.v.e.p) then the local resolvent of every vector $x\in\x$ is the set
of complex numbers $\lambda$ such that there exists an vector valued analytic function $F$ defined in a open neighborhood $V$ of
$\lambda$ which verifies $(T-\mu I)F(\mu )=x~~\mbox{for every }\mu\in V.$\\
{\bf Example of an operator without s.v.e.p. }Let $T$ be the adjoint of the unilateral weighted shift operator with the weight
$\omega_n=1$ for every $n\in\N$. i.e:

\begin{displaymath}
Te_n=\left\{
\begin{array}{lll}
0&\textrm{\mbox{if }$n=0$}\\
\\
e_n&\textrm{\mbox{if }$n>0$}\\
\end{array}
\right.
\end{displaymath}
Clearly, we have $\|T\|=1$ and \sst =$\big{\{}\lambda\in\C~:~|\lambda |\leq 1\big{\}}$ then \rrt =$\big{\{}\lambda\in\C~:~|\lambda |>
1\big{\}}.$
Therefore, for every $\lambda\in$\rrt~we have $\|\frac{T}{\lambda}\|<1$ and so,
\begin{eqnarray*}
(T-\lambda I)^{-1}&=&\frac{1}{\lambda}(\frac{T}{\lambda}- I)^{-1}\\
&=&\frac{1}{\lambda}\sum\limits_{n=0}^{+\infty}(\frac{T}{\lambda})^n\\
&=&-\sum\limits_{n=0}^{+\infty}\frac{T^n}{\lambda^{n+1}}.
\end{eqnarray*}
Hence for every $\lambda\in\rho (T),$
\begin{eqnarray*}
R_{_{T,e_0}}&=&(T-\lambda I)^{-1}e_0\\
&=&-\sum\limits_{n=0}^{+\infty}\frac{T^ne_0}{\lambda^{n+1}}\\
&=&-\frac{e_0}{\lambda} ~~\mbox{because }T^ne_0=0\mbox{ for every }n\geq 0.
\end{eqnarray*}
We will now exhibit two different analytic extensions of $R_{_{T,e_0}}$. First let $F(\lambda )=-\frac{e_0}{\lambda}$ for
$\lambda\not= 0;$ then $(T-\lambda I)F(\lambda )=e_0$ for every $\lambda\not= 0$. And so, $F$ is clearly an analytic extension of
$R_{_{T,e_0}}$. Second, set
\begin{displaymath}
G(\lambda )=\left\{
\begin{array}{lll}
-\frac{e_0}{\lambda}&\textrm{\mbox{if }$|\lambda |>1$}\\
\\
\sum\limits_{n=0}^{+\infty}\lambda^ne_{n+1}&\textrm{\mbox{if }$|\lambda |<1$}.\\
\end{array}
\right.
\end{displaymath}
Clearly, that $G$ has the right form for $\lambda\in\C$, $|\lambda |>1$. For $\lambda\in\C$, $|\lambda |<1$ we have
\begin{eqnarray*}
(T-\lambda I)G(\lambda )&=&(T-\lambda I)\sum\limits_{n=0}^{+\infty}\lambda^ne_{n+1}\\
&=&\sum\limits_{n=0}^{+\infty}\lambda^nTe_{n+1}-\sum\limits_{n=0}^{+\infty}\lambda^{n+1}e_{n+1}\\
&=&\sum\limits_{n=0}^{+\infty}\lambda^ne_n-\sum\limits_{n=0}^{+\infty}\lambda^{n+1}e_{n+1}\\
&=&e_0.
\end{eqnarray*}
Thus $G$ is also an analytic extension of $R_{_{T,e_0}}.$ However the two analytic extensions $F$ and $G$ do not coincide on
$\big{\{}\lambda\in\C~:~0<|\lambda |<1\big{\}}.$\\
{\bf Note. }The operator $T$ considered in this example is the adjoint of the of the unilateral weighted shift operator with the
weight $\omega_n=1$ for every $n\in\N$ which is an isometric non unitary operator. More than that one can prove that the adjoint
operator of every isometric non unitary operator on Hilbert space do not has the s.v.e.p (see \cite{8}).\\
{\bf Example of operators having the s.v.e.p.} Every operator $T\in\lx$ which has empty interior of its point spectrum has the
s.v.e.p. Let $F$ be a vector valued analytic function such that
$$(T-\lambda I)F(\lambda )=0~~\mbox{for every }\lambda\in D(F).$$
Then $$TF(\lambda )=\lambda F(\lambda )~~\mbox{for every }\lambda\in D(F).$$
And so, $F$ must be identically zero otherwise there is $\lambda\in D(F)$ and $r>0$ such that
$B_r=\big{\{}\mu\in\C~:~|\lambda -\mu |<r\big{\}}\subset D(F)$ and
$F(\mu )\not= 0$ for every $\mu\in B_r.$ Thus, $B_r\subset \sigma_p(T)$. Contradiction.
Hence, $F\equiv 0$ and so, $T$ has the s.v.e.p.\\
{\bf Note. }From Theorem 1.3.12 and Remark 1.3.13 follow that every unilateral weighted shift operator has the s.v.e.p.\\
We will see that all operators in which we are interested have the s.v.e.p. So, from now on we will be dealing with operators which
have only single-valued extension property (s.v.e.p).
\begin{prop}Let $T\in\lx$ be an operator which has the s.v.e.p and let $x,~y$ be two vectors of \x.~Then:\\
{\rm (i)}~~$\lsx (\alpha x)=\lsx (x)$ for every $\alpha\not= 0$.\\
{\rm (ii)}~$\lsx (x+y)\subset\lsx (x)\cup\lsx (y)$.\\
{\rm (iii)}$\lsx (x)$ is empty if and only if $x=0$.\\
{\rm (iv)}~$\lsx \big{(}\widetilde{x}(\lambda )\big{)}=\lsx (x)$ for every $\lambda\in\lrx (x)$ fixed point.\\
{\rm (v)}~~For every operator $S\in\lx$ which commutes with $T$, $\lsx (Sx)\subset\lsx (x)$.
\end{prop}
{\bf Proof. }(ii)~~The result follows immediately since
$$(T-\lambda I)\big{(}\widetilde{x}(\lambda )+\widetilde{y}(\lambda )\big{)}=x+y~~\mbox{for every }\lambda\in\lrx (x)\cap\lrx (y).$$
(iii)~If $x=0$ then there noting to prove. Now, suppose that $\lsx (x)=\emptyset$ then the maximal single valued extension
$\widetilde{x}$ is an entire function which coincide with $R_{_{T,x}}$ on \rrt.~
Since $\big{\{}\lambda\in\C ~:~|\lambda |>\|T\|\big{\}}\subset\rho (T)$ and
$\widetilde{x}(\lambda )=R_{_{T,x}}(\lambda )=R_{_T}(\lambda )x$ for $|\lambda |>\|T\|$ then from the Theorem 1.2.1 it follows that
$\lim\limits_{|\lambda |\to +\infty}\widetilde{x}(\lambda )=0.$ Therefore by Liouville's Theorem $$\widetilde{x}\equiv 0.$$
In particular, $R_{_{T,x}}\equiv 0$ on \rrt.~ Hence $x=0.$\\
(iv)~~Fix $\lambda\in\lrx (x)$ and let $y=\widetilde{x}(\lambda ).$ Let we first prove that $\lsx (x)\subset\lsx (y).$ We have
$$(T-\mu I)\widetilde{y}(\mu )=y \mbox{ for every }\mu\in\lrx (y)$$
then for every $\mu\in\lrx (y)$ we have
\begin{eqnarray*}
(T-\lambda I)(T-\mu I)\tilde{y}(\mu )&=&(T-\mu I)(T-\lambda I)\tilde{y}(\mu )\\
&=&(T-\lambda I)y\\
&=&(T-\lambda I)\widetilde{x}(\lambda )~~\mbox{because }y=\widetilde{x}(\lambda )\\
&=&x
\end{eqnarray*}
Since $(T-\lambda I)\widetilde{y}$ is a vector valued analytic function on $\lrx (y)$ then it is analytic extension of $R_{_{T,x}}$.
Therefore, $\lrx (y)\subset\lrx (x)$. Thus $\lsx (x)\subset\lsx (y)=\lsx \big{(}\widetilde{x}(\lambda )\big{)}.$\\
Conversely, let use to prove that $\lsx \big{(}\widetilde{x}(\lambda )\big{)}\subset\lsx (x).$ Clearly that the vector valued
function $F$ defined in $\lrx (x)$ by
\begin{displaymath}
F(\lambda )=\left\{
\begin{array}{lll}
\frac{\widetilde{x}(\mu )-\widetilde{x}(\lambda )}{\mu -\lambda }&\textrm{\mbox{if }$\mu\not=\lambda$}\\
\\
\widetilde{x}^\prime (\lambda )&\textrm{\mbox{if }$\mu =\lambda$}
\end{array}
\right.
\end{displaymath}
is an analytic function on $\lrx (x)$ and for $\mu\not=\lambda$ we have
\begin{eqnarray*}
(T-\mu I)F(\mu )&=&(T-\mu I)\bigg{[}\frac{\widetilde{x}(\mu )-\widetilde{x}(\lambda )}{\mu -\lambda }\bigg{]}\\
&=&\frac{1}{\mu -\lambda}\bigg{[}(T-\mu I)\widetilde{x}(\mu )-(T-\mu I)\widetilde{x}(\lambda )\bigg{]}\\
&=&\frac{1}{\mu -\lambda}\bigg{[}x-(T-\lambda I)\widetilde{x}(\lambda )+(\mu -\lambda )\widetilde{x}(\lambda )\bigg{]}\\
&=&\frac{1}{\mu -\lambda}\bigg{[}x-x+(\mu -\lambda )\widetilde{x}(\lambda )\bigg{]}\\
&=&\widetilde{x}(\lambda )\\
&=&y
\end{eqnarray*}
This equality holds also for $\mu =\lambda$ by making $\mu\to\lambda$ in the last equality, and so
$$(T-\mu I)F(\mu )=\widetilde{x}(\lambda )~~\mbox{for every }\lambda\in\lrx (x).$$
Hence, $\lrx (x)\subset\lrx \big{(}\widetilde{x}(\lambda )\big{)}.$
Thus
$\lsx \big{(}\widetilde{x}(\lambda )\big{)}\subset\lsx (x).$
And so,
$\lsx (x)=\lsx \big{(}\widetilde{x}(\lambda )\big{)}~~\mbox{for every }\lambda\in\lrx (x).$\\
(v)~~~We have
$(T-\lambda I)\widetilde{x}(\lambda)=x~~\mbox{for every }\lambda\in\lrx (x)$
then $S(T-\lambda I)\widetilde{x}(\lambda)=Sx~~\mbox{for every }\lambda\in\lrx (x)$
and so,
$$(T-\lambda I)S\widetilde{x}(\lambda)=Sx~~\mbox{for every }\lambda\in\lrx (x).$$
Hence the vector valued function $S\widetilde{x}$ is an analytic extension of $R_{_{T,Sx}}$. Therefore, $\lrx (x)\subset\lrx (Sx).$  
Thus, $\lsx (Sx)\subset\lsx (x)$.$\hfill\blacksquare$
\begin{thm}
For every operator $T\in\lx$ which has the s.v.e.p
$$\sigma (T)=\bigcup\limits_{x\in{\cal X}}\lsx (x).$$
\end{thm}
{\bf Proof. }Assume that there is $\lambda\in\sigma (T)$ such that $\lambda\not\in\lsx (x)$ for every $x\in\x.$
Then $\lambda\in\bigcap\limits_{x\in{\cal X}}\lrx (x)$ and so, $(T-\lambda I)\tilde{x}(\lambda )=x$ for every $x\in\x,$ in
particular the
operator $T-\lambda I$ is surjective. Therefore it follows from the open mapping Theorem (see \cite{6} and \cite{7}) that $T-\lambda I$ must be not injective
since $\lambda\in\sigma (T)$. So, there is $x_1\not= 0\in\x$ such that $(T-\lambda I)x_1=0.$ Since $T-\lambda I$ is bounded surjective
operator then it follows from the Open Mapping Theorem that there is a constant $c>0$ such that for every $x\in$\x~there exists
$y\in$\x~ with
$$\|y\|<c\|x\|~~\mbox{and}~~(T-\lambda I)y=x.$$
Therefore, there is a sequence $(x_n)_n$ of elements of \x~ such that
$$(T-\lambda I)x_{n+1}=x_n~~\mbox{and }\|x_{n+1}\|\leq c\|x_n\|~~\mbox{for every }n\geq 1 .$$
And so, the function $F$ defined on the open disc $D=\big{\{}\mu\in\C~:~|\mu-\lambda |<\frac{1}{c}\big{\}}$ by 
$$F(\mu )=\sum\limits_{n=0}^{\+\infty}x_{n+1}(\mu -\lambda )^n$$
is a non identically zero analytic function on $D$ and for every $\mu\in D$ we have
\begin{eqnarray*}
(T-\mu I)F(\mu )&=&(T-\lambda )F(\mu )+(\lambda -\mu )F(\mu )\\
&=&\sum\limits_{n=0}^{+\infty}(T-\lambda )\bigg{(}x_{n+1}(\mu -\lambda )^n\bigg{)}-(\lambda -\mu
)\sum\limits_{n=0}^{+\infty}x_{n+1}(\mu -\lambda )^n\\
&=&\sum\limits_{n=1}^{+\infty}x_n(\mu -\lambda )^n-\sum\limits_{n=0}^{+\infty}x_{n+1}(\mu -\lambda )^{n+1}\\
&=&0
\end{eqnarray*}
Contradiction with the assumption that $T$ has the s.v.e.p, and so the desired result holds.$\hfill\blacksquare$
\begin{lem}Let $T\in$\lx~be a surjective bounded operator which has the s.v.e.p. Then $T$ is invertible in
\lx.~Therefore,
$$\sigma (T)=\big{\{}\lambda\in\C~:~T-\lambda I~\mbox{is not onto}\big{\}}.$$
\end{lem}
{\bf Proof. }Suppose that $T$ is not invertible in \lx~then it is not injective. So, there is an element $x_1\not= 0\in {\cal X}$
such that $Tx_1=0$. And so, by the same argument of the proof of Theorem 2.1.2 it follows that $T$ has not the s.v.e.p.$\hfill\blacksquare$
\begin{thm}Let $T\in\lh$ be an bounded operator on a Hilbert space \h.~If its adjoint $T^*$ has s.e.v.p then
$$\sigma (T)=\sigma_{ap}(T).$$
\end{thm}
{\bf Proof. }It is known that $\sigma_{ap}(T)\subset\sigma (T).$ Conversely, let $\lambda\not\in\sigma_{ap}(T)$ then
$T-\lambda I$ is bounded below. So, it follows from Theorem 1.1.8 that $T^*-\overline{\lambda} I$ is surjective. Hence by the
Lemma 2.1.3 $\overline{\lambda}\not\in\sigma (T^*)$. Thus $\lambda\not\in\sigma (T)$ since $\sigma (T)=\overline{\sigma (T^*)}$.
This complete the proof of the Theorem.$\hfill\blacksquare$
\section{Operators with the Dunford's Condition C (DCC)}
\hspace*{1.00cm}Let $T\in\lx$ be an operator which has s.v.e.p. So,
from the properties (i) and (ii) of Proposition 2.1.1 it follows that for every closed subset $F$ of \C~the following
subset of \x ~given by
$${\cal X}_{_T}(F)=\big{\{}x\in{\cal X}~:~\lsx (x)\subset F\big{\}}$$
is a linear subspace of \x~not necessarily closed as we will see in the example below.
The operator $T$ is said to satisfy the {\bf Dunford's Condition C (DCC)} if for every closed subset $F$ of \C~the linear subspace
$\xtf$ is closed. We will see that all the operators in which we are interested have the Dunford's Condition C (DCC).
Also, from the property (v) of the same proposition we see that the linear subspace ${\cal X}_{_T}(F)$ is invariant linear subspace
for every operator $S\in\lx$ which commutes  with $T$. i.e:
$$S\bigg{(}{\cal X}_{_T}(F)\bigg{)}\subset {\cal X}_{_T}(F)~~~\mbox{for every }S\in{\cal L(X)}, ~\mbox{with } ST=TS.$$
If the linear subspace $\xtf$ is a closed and not trivial (i.e: $\xtf\not= \{0\}$ and $\xtf\not=\x$) then it is a non trivial
hyperinvariant linear subspace for $T$ and this is one method to construct the non trivial hyperinvariant subspaces for a given
operator in \lx~. (More information about the existence of hyperinvariant subspaces of operators can be found in \cite{1}, \cite{3} and \cite{14}).\\
Since for every $x\in\x$ we have $\lsx (x)\subset\sigma (T)$ then obviously we have $\xtf =\xxt \big{(}F\cap\sigma (T)\big{)}$,
hence $F$ can be supposed to be a closed subset of \sst .\\
Before proceeding further it would be well to point out that the Dunford's Condition C does not hold for every operator which has the
single valued extension property. Let us consider an example. Let $(e_n)_n$ be an orthonormal basis for an Hilbert space \h,~and let
$T$ be the unilateral weighted shift operator given by
\begin{displaymath}
Te_n=\left\{
\begin{array}{lll}
e_{n+1}&\textrm{\mbox{if }$n$\mbox{ is not square integer,}}\\
\\
0&\textrm{\mbox{if }$n$\mbox{ is square integer.}}\\ 
\end{array}
\right.
\end{displaymath}
It is easy to see that $\|T^n\|=1$ for every $n\geq 1$, then spectral radius
$r(T)$ of $T$ is $1$. On the other hand, $T$ has the s.v.e.p property since $\sigma_p(T)=\{0\}$ (see Remark 1.3.13)
Now let us suppose that ${\cal X}_{_T} (\{0\})$ is closed. For every non-negative integer $n$ there is 
an integer $k$ such that $k^2\leq n<(k+1)^2$ and so,
\begin{itemize}
\item If $n=k^2$ then $$(T-\lambda I)\bigg{(}-\frac{e_{n}}{\lambda}\bigg{)}=e_{n}~~~\mbox{for }\lambda\not= 0.$$ 
\item If $n=k^2+s$ for some $0<s<2k+1$ then 
$$(T-\lambda I)\bigg{(}-\frac{e_{n}}{\lambda}-\frac{e_{n+1}}{\lambda^2}-...-\frac{e_{(k+1)^2}}{\lambda^{2(k+1)-s}}\bigg{)}=e_{n}~~~\mbox{for }\lambda\not= 0.$$
\end{itemize}
Then for every integer $n\geq 0$, $\lsx (e_n)=\{0\}$. And so, ${\cal X}_{_T}(\{0\})=$\h.~ Hence it follows from Theorem 2.1.2
that
$\sigma (T)=\{0\}$. This is impossible since the spectral radius $r(T)$ of $T$ is $1$.
\begin{thm}For every closed subset $F$ of $\sigma (T)$ such that the linear subspace $\xtf$ is closed,
$$\sigma \big{(}T_{|\xtf}\big{)}\subset F.$$
\end{thm}
{\bf Proof.}
Let $\lambda\not\in F$. We will prove that $\lambda\in\rho\big{(}T_{|\xtf}\big{)}$, that is $T_{|\xtf}-\lambda I$ is invertible in
${\cal L(\xtf )}.$
Since for every $x\in\xtf$ we have $\C\backslash F\subset\lrx (x).$ Therefore $\widetilde{x}(\lambda )$
makes sense and from the equality (iv) of Proposition 2.1.1 it follows that $\lrx (x)=\lrx (\widetilde{x}(\lambda )).$ Thus
$\widetilde{x}(\lambda )\in\xtf.$ Let $A$ be the map from $\xtf$ to it self defined by $Ax=\widetilde{x}(\lambda ).$ It is evident
that $A$ is linear map. We will show by using the Closed Graph Theorem (see \cite{6} and \cite{7}) that the linear map $A$ is bounded. Let
$(x_n)_n$ be a sequence of elements of $\xtf$ such that
$\lim\limits_nx_n=x\in{\cal X}$ and $\lim\limits_nAx_n=y\in{\cal X}$. Since $\xtf$ is closed linear subspace of \x~then
$x\in\xtf$ and $y=\lim\limits_nAx_n=\lim\limits_n\widetilde{x_n}(\lambda )\in\xtf.$ We have
$$(T-\lambda I)Ax_n=(T-\lambda I)\widetilde{x_n}(\lambda )=x_n,$$
from which it follows by the continuity of the linear operator $(T-\lambda I)$ that $(T-\lambda I)y=x.$ On the other hand,
$(T-\lambda I)\widetilde{x}(\lambda )=x,$
therefore,
$(T-\lambda I)\bigg{(}\widetilde{x}(\lambda )-y\bigg{)}=0.$
Since $\xtf$ is a linear subspace and $y,~\widetilde{x}(\lambda )\in\xtf$ then $z=\widetilde{x}(\lambda )-y\in\xtf.$ Hence
$$\sigma_{_T}(z)\subset F.\leqno{(*)}$$
Let $G$ be the vector valued analytic function defined on $\C\backslash\{\lambda\}$ by $G(\mu )=\frac{1}{\lambda -\mu }z$.
Then for every $\mu\in\C\backslash\{\lambda\}$ we have
\begin{eqnarray*}
(T-\mu I)G(\mu )&=&\frac{1}{\lambda -\mu }(T-\mu I)z\\
&=&\frac{1}{\lambda -\mu }\bigg{(}(T-\lambda I)z+(\lambda -\mu )z\bigg{)}\\
&=&\frac{1}{\lambda -\mu }\bigg{(}0+(\lambda -\mu )z\bigg{)}\\
&=&z.
\end{eqnarray*}
Hence $\C\backslash\{\lambda\}\subset\lrx (z).$ Thus
$$\lsx (z)\subset \{\lambda\}.\leqno{(**)}$$
Since $\lambda\not\in F$ then from $(*)$ and $(**)$ it follows that $\lsx (z)\subset F\cap \{\lambda\}=\emptyset .$ Hence
$\lsx (z)=\emptyset$ and so, $z=0$ by the property (iii) of Proposition 2.1.1. Therefore,
$$Ax=\widetilde{x}(\lambda )=y.$$
And so, by the Closed Graph Theorem the linear operator $A$ is bounded on $\xtf$.\\
Now, let us prove that the operator $A$ is exactly the inverse of $T_{|\xtf}-\lambda I$. For every $x\in\xtf$ we have
\begin{eqnarray*}
(T_{|\xtf}-\lambda I)Ax&=&(T-\lambda I)\widetilde{x}(\lambda )\\
&=&x.
\end{eqnarray*}
On the other hand for every $x\in\xtf$, we have $\widetilde{(T-\lambda I)x}(\lambda )=(T-\lambda I)\widetilde{x}(\lambda )$, for which,
according to the definition of $A$, it follows that
\begin{eqnarray*}
A(T_{|\xtf}-\lambda I)x&=&\widetilde{(T-\lambda I)x}(\lambda )\\
&=&(T-\lambda I)\widetilde{x}(\lambda )\\
&=&x.
\end{eqnarray*}
Therefore $\lambda\in\rho\big{(}T_{|\xtf}\big{)}$.$\hfill\blacksquare$\\
\\
\hspace*{1.00cm}We shall see an application of this Theorem. The operator $T$ is said to be {\bf cyclic} if there is a vector $x\in$\x~such that
the linear subspace generated by $\{T^nx~:~n\in\N\}$ is dense in \x.~The vector $x$ is called a cyclic vector for $T$. If $T$ is an
injective unilateral weighted shift with weights $(\omega_n)_n$ that is $Te_n=\omega_ne_{n+1},~n\geq 0$ where $(e_n)_n$ is a
orthonormal basis of \h,~then $T$ is cyclic with cyclic vector $e_0$ since $T^ne_0=\omega_0...\omega_{n-1}e_n$ for every $n\geq 1$.
(For more information about cyclic shifts see \cite{25} and \cite{255}).
\begin{prop}Suppose that the operator $T$ is cyclic with a cyclic vector $x\in$\x.~If the operator $T$ satisfies DCC then
$\lsx (x)=\sigma (T).$
\end{prop}
{\bf Proof. }Let $F=\lsx (x)$. Then $x\in\xtf$. Since $\xtf$ is invariant subspace for $T$ then
$T^nx\in\xtf$ for every $n\in\N.$ Hence, the linear subspace generated by $\{T^nx~:~n\in\N\}$ is contained in $\xtf$. Since the
operator $T$ satisfies DCC property then $\xtf$ is closed subspace and so it follows from the density in \x~of the linear subspace
generated by $\{T^nx~:~n\in\N\}$ that \x $=\xtf$. And so, from Theorem 2.2.1 it follows that
$$\sigma (T)=\sigma (T_{_{|\xtf}})\subset F=\lsx (x).$$
Therefore, $\sigma (T)=\lsx (x)$.$\hfill\blacksquare$
\begin{rem} In this proposition we did not use the DDC property of the operator $T$. We used only the fact that $\xtf$ is closed
linear subspace where $F=\lsx (x)$.
\end{rem}

\chapter{Subnormal and Hyponormal Operators}
\hspace*{1.00cm}Several classes of Hilbert spaces operators are defined around the notion of normal operator. In 1950, Paul R. Halmos, motivated by
the successful development of the theory of normal operators, introduced the notions of subnormality and hyponormality for bounded
Hilbert space operators, in an attempt to extend the basic facts of the spectral theory of normal operators. Those classes of operators
have been the subject of much investigation during the last fifty years and many important developments in the operator theory
have dealt with them e.g: S.Brown's proof of the existence of non trivial subspaces, J.Conway and R.Olin's construction of the functional calculus
and J.Thomson's description of the spectral picture in the cyclic case for subnormal operators, etc. Our goal in this chapter is
to study the structure of those classes of operators and to collect sufficient material to study the local spectra of cyclic hyponormal
operators in the next chapter.\\  
\section{Subnormal and Hyponormal Operators}
\hspace*{1.00cm}In this section, we recall first the general and fundamental properties of normal bounded Hilbert space operators.
Let \lh~be the algebra of all linear bounded operators on a Hilbert space \h.~An operator $T\in$\lh~ is said to be {\bf normal} if
it commutes with its adjoint $T^*$ i.e: $T^*T=TT^*.$\\
We first note that for every operator $T\in$\lh,
\begin{equation}
\big{\langle}\big{(}T^*T-TT^*\big{)}x~,~x\big{\rangle} =\|Tx\|^2-\|T^*x\|^2~~\mbox{for every }x\in{\cal H}.
\end{equation}
\begin{prop}An operator $T\in$\lh~ is normal if and only if $\|Tx\|=\|T^*x\|$ for every $x\in$\h.
\end{prop}
{\bf Proof. }Since the zero operator is the unique operator $S\in$\lh~which has the property that
$$\big{\langle} Sx~,~x\big{\rangle} =0  ~~~\mbox{for every }x\in{\cal H},$$
then the proof follows immediately from the identity (3.1).$\hfill\blacksquare$\\
\hspace*{1.00cm}Using Liouville's Theorem, Fugled Putman in 1950 give a nice characterization of normal operators.
\begin{thm}[Fuglede] An operator $T\in\lh$ is normal if and only if $T^*S=ST^*$ for every operator $S\in\lh$ such that $TS=ST$.
\end{thm}
\hspace*{1.00cm}In the proof of this theorem, we shall need the following:\\
Let $T$ be an operator in \lh.~For every positive integer $n$, let $T_n=I+\frac{T}{1!}+\frac{T^2}{2!}+\ldots+\frac{T^n}{n!}$.
$T_n$ is a bounded operator in \lh;~and the sequence $(T_n)_n$ converges in \lh~to an invertible operator denoted by $e^T$ and its
inverse is $e^{-T}$ i.e:
\begin{equation}
\bigg{(}e^T\bigg{)}^{-1}=e^{-T};
\end{equation}
in addition, if $S$ is an operator in \lh~which commutes with $T$ then:
\begin{equation}
e^TS=Se^T,
\end{equation}
and
\begin{equation}
e^{T+S}=e^Te^S=e^Se^T.
\end{equation}
{\bf Proof of Theorem 3.1.2. }Let $S$ be an operator in \lh~such that $TS=ST.$ Then it follow from (3.2) and (3.3) that
$S=e^{-i\overline{z}T}Se^{i\overline{z}T}~~\mbox{for every }z\in\C.$ And so, for every $z\in\C$,
\begin{eqnarray*}
e^{-izT^*}Se^{izT^*}&=&e^{-izT^*}e^{-i\overline{z}T}Se^{i\overline{z}T}e^{izT^*}\\
&=&e^{-i\big{(}zT^*+\overline{z}T\big{)}}Se^{i\big{(}zT^*+\overline{z}T\big{)}}
\end{eqnarray*}
Since the operator $e^{i\big{(}zT^*+\overline{z}T\big{)}}$ is unitary and its adjoint is
$e^{-i\big{(}zT^*+\overline{z}T\big{)}}$ then
$$e^{-izT^*}Se^{izT^*}=\bigg{(}e^{i\big{(}zT^*+\overline{z}T\big{)}}\bigg{)}^*Se^{i\big{(}zT^*+\overline{z}T\big{)}};$$
therefore the following operator valued function defined on \C~by $\phi (z)=e^{-izT^*}Se^{izT^*}$ is bounded analytic
function, by Liouville's Theorem it is a constant function; in particular its derivative $\phi^\prime$ is zero. So,
\begin{eqnarray*}
\phi^\prime (z)&=&-iT^*e^{-izT^*}Se^{izT^*}+e^{-izT^*}S(iT^*)e^{izT^*}\\
&=&-iT^*\phi (z)+e^{-izT^*}Se^{izT^*}(iT^*)~~~~\mbox{by }(1.5)\\
&=&-iT^*\phi (z)+\phi (z)(iT^*)\\
&=&0.
\end{eqnarray*}
Hence,$$T^*\phi (z)=\phi (z)T^*~~~~\mbox{for every }z\in\C.$$ Thus $T^*S=ST^*$ because $\phi (0)=S$.$\hfill\blacksquare$\\
\\
\hspace*{1.00cm}Let $T$ be the unweighted bilateral shift i.e: $Te_n=e_{n+1}$ for every $n\in\Z$ where $(e_n)_n$ is an orthonormal
basis for \h.~Then $T$ is unitary operator since by Proposition 1.3.2, we have $T^*e_n=e_{n-1}$ for every $n\in\Z$. So,
in particular, it is normal operator. In fact, one can see that a bilateral weighted shift operator $T$ with a
nonnegative weights $(\omega_n)_{n\in\Z}$ is normal if and only if the weights $(\omega_n)_{n\in\Z}$ are constant; and a non
zero unilateral shift is never normal.\\
We note that the restriction of a normal operator $T\in$\lh~ on a proper closed invariant subspace is not in general normal
operator. Let $S$ be the restriction of the unweighted bilateral shift $T$ on the closed linear subspace $K$ generated by
$\{e_n~:~n\geq 0\}$ which is proper closed invariant subspace for $T$. So, the restriction $S$ is exactly the unweighted
unilateral shift on the Hilbert space ${\cal K}$ which is not normal but it has a normal extension. Therefore  it is natural to
introduce  and study the theory of subnormal operators which constitute a considerably more useful and deeper generalization of
the theory of the normal operators. An operator $T\in\lh$ is called a {\bf subnormal operator} if it has a normal extension i.e: if
there exists a normal operator $S$ on a Hilbert space ${\cal K}$ such that $\h$ is a closed invariant subspace for $S$ and the
restriction of $S$ to \h~ coincides with $T$. A normal extension $S$ on a Hilbert space ${\cal K}$ of an subnormal operator
$T\in\lh$ is called {\bf minimal normal extension} if there is no closed invariant subspace for $S$ on which the restriction
of $S$ is normal operator; a such extension always exists by Zorn's Lemma and it is a unique up to an invertible isometry.
\begin{lem}Let $T\in\lh$ be a subnormal operator. A normal extension $S$ of $T$ on an Hilbert space ${\cal K}$ is a minimal normal
extension if and only if ${\cal K}$ coincide with the closure of the linear subspace generated by $\{S^{*^n}x~:~n\in\N\}$.
\end{lem}
{\bf Proof. }It suffices to observe that the closed linear subspace ${\cal K}_0$ generated by $\{S^{*^n}x~:~n\in\N\}$ is containing $\h$, 
invariant by $S$ and $S^*$ and the restriction of $S$ over ${\cal K}_0$ is normal.$\hfill\blacksquare$
\begin{lem}Let $S$ be a normal extension on a Hilbert space ${\cal K}$ of an subnormal operator $T\in\lh$. Then for every finite
family of elements $x_1,x_2,...,x_n\in\h$ we have,
$$\|\sum\limits_{i=1}^nS^{*^i}x_i\|=\|\sum\limits_{i=1}^nT^{*^i}x_i\|.$$
\end{lem}
{\bf Proof. }We have
\begin{eqnarray*}
\|\sum\limits_{i=1}^nS^{*^i}x_i\|^2
&=&\bigg{\langle}\sum\limits_{i=1}^nS^{*^i}x_i~,~\sum\limits_{j=1}^nS^{*^j}x_j\bigg{\rangle}\\
&=&\sum\limits_{1\leq i,j\leq n}\big{\langle}S^{*^i}x_i~,~S^{*^j}x_j\big{\rangle}\\
&=&\sum\limits_{1\leq i,j\leq n}\big{\langle}S^jx_i~,~S^ix_j\big{\rangle}\\
&=&\sum\limits_{1\leq i,j\leq n}\big{\langle}T^jx_i~,~T^ix_j\big{\rangle}~~~~~~~\mbox{because }S_{|\h}=T\\
&=&\|\sum\limits_{i=1}^nT^{*^i}x_i\|^2.
\end{eqnarray*}
$\hfill\blacksquare$
\begin{thm}If $S_1$ and $S_2$ are two minimal normal extensions, on Hilbert spaces respectively ${\cal K}_1$ and ${\cal K}_2$, of
a subnormal operator $T\in\h$ then there exists an invertible isometry $U$ from ${\cal K}_1$ onto ${\cal K}_2$ such that  
$US_1=S_2U$ and $U(x)=x$ for every $x\in\h$.
\end{thm}
{\bf Proof. }Let $M_1$ and $M_2$ be the linear subspaces respectively of ${\cal K}_1$ and ${\cal K}_2$ generated respectively by the sets
$\{S_1^{*^i}x~:~i\in\N~\mbox{and }x\in\h\}$ and $\{S_2^{*^i}x~:~i\in\N~\mbox{and }x\in\h\}$. It follows from the Lemma 3.1.3 that
${\cal K}_1=\overline{M_1}$ and ${\cal K}_2=\overline{M_2}$.
On the other, it follows from the Lemma 3.1.4 that for every finite family of elements $x_1,x_2,...,x_n\in\h$
$$\|\sum\limits_{i}S_1^{*^i}x_i\|=\|\sum\limits_{i}S_2^{*^i}x_i\|.$$
Hence the correspondence defined from $M_1$ onto $M_2$ by $\sum\limits_{i}S_1^{*^i}x_\longmapsto \sum\limits_{i}S_2^{*^i}x_i$ is an
isometry, we denote it by $U$. Therefore $U$ has a unique isometric extension that maps ${\cal K}_1$ onto ${\cal K}_2$ and it is   
the identity on $\h$. To prove that $US_1=S_2U$, it suffices to verify that $US_1$ agrees with $S_2U$ on $M_1$ and this is implied
by
\begin{eqnarray*}
US_1\bigg{(}\sum\limits_{i}S_1^{*^i}x_i\bigg{)}&=&U\bigg{(}\sum\limits_{i}S_1S_1^{*^i}x_i\bigg{)}
=U\bigg{(}\sum\limits_{i}S_1^{*^i}S_1x_i\bigg{)}\\
&=&U\bigg{(}\sum\limits_{i}S_1^{*^i}Tx_i\bigg{)}
=\bigg{(}\sum\limits_{i}S_2^{*^i}Tx_i\bigg{)}\\
&=&\bigg{(}\sum\limits_{i}S_2^{*^i}S_2x_i\bigg{)}
=S_2\bigg{(}\sum\limits_{i}S_2^{*^i}x_i\bigg{)}\\
&=&S_2U\bigg{(}\sum\limits_{i}S_1^{*^i}x_i\bigg{)}.
\end{eqnarray*}
This complete the proof of the Theorem.$\hfill\blacksquare$
\begin{thm}If $T\in\lh$ is a subnormal operator and $S$ is its minimal normal extension on a Hilbert space ${\cal K}$ then $\sigma
(S)\subset\sigma (T)$. In particular, the minimal normal extension of an invertible subnormal operator is invertible.
\end{thm}
{\bf Proof. }To prove that $\sigma (S)\subset\sigma (T)$ it suffices to prove that $S-\lambda I$ is invertible for every
$\lambda\in\rho (T)$. On the other hand, $S-\lambda I$ is a minimal normal extension of the operator $T-\lambda I$ so, the assertion
reduces to prove that if $T$ is invertible then its minimal normal extension $S$ is invertible.
It is clear that the closure of the range of $S$ is a closed invariant subspace $M$ for $S$ and for its adjoint $S^*$. Therefore the
restriction of $S$ on $M$ is normal. On the other hand, $\h =T\h=S\h$ is contained in the range of $S$. Hence it follows from the
minimality of $S$ that $M={\cal K}$. So, $S$ has dense range, therefore it suffices to prove that $S$ is bounded from below.
Let $L$ be the linear subspace generated by the set $\{S^{*^i}x~:~i\N~\mbox{and }x\in\h\}$ then by the Lemma 3.1.3 it
follows that $\overline{L}={\cal K}$. For every finite sum $\sum\limits_{i}S^{*^i}x_i\in L$ we have
\begin{eqnarray*}
\|S\bigg{(}\sum\limits_{i}S^{*^i}x_i\bigg{)}\|&=&\|S^*\bigg{(}\sum\limits_{i}S^{*^i}x_i\bigg{)}\|~~~~\mbox{by proposition
}3.1.1\\ 
&=&\|\sum\limits_{i}S^{*^{i+1}}x_i\|\\
&=&\|\sum\limits_{i}T^{*^{i+1}}x_i\|~~~~~~~~~~\mbox{Lemma (3.1.4)}\\
&=&\|T^*\bigg{(}\sum\limits_{i}T^{*^{i}}x_i\bigg{)}\|\\
&\geq &\frac{1}{\|T^{*^{-1}}\|}\|\sum\limits_{i}T^{*^i}x_i\|=\frac{1}{\|T^{*{-1}}\|}\|\sum\limits_{i}S^{*^i}x_i\|~~~~\mbox{Lemma (3.1.4)}.
\end{eqnarray*}
Thus $\frac{1}{\|T^{*^{-1}}\|}\|x\|\leq \|Sx\|~~~\mbox{for every }x\in {\cal K};$ and the desired result follows.$\hfill\blacksquare$
\begin{thm}For every subnormal operator $T\in\lh$, $\|T^*x\|\leq \|Tx\|$ for every $x\in\h$.
\end{thm}
{\bf Proof. }First let us prove that $T^*x=PS^*x$ for every $x\in\h$ where $S$ a minimal normal extension on a Hilbert space
${\cal K}$ and $P$ is the linear projection from ${\cal K}=\h\oplus \h^\perp$ onto $\h$. For every $x,~y\in\h$ we have
\begin{eqnarray*}
\langle T^*x~,~y\rangle&=&\langle x~,~Ty\rangle =\langle x~,~Sy\rangle ~~~~~~~~~\mbox{because }S_{|\h}=T\\
&=& \langle S^*x~,~y\rangle =\langle S^*x~,~Py\rangle~~~~~\mbox{because }P_{|\h}=I\\
&=&\langle P^*S^*x~,~y\rangle =\langle PS^*x~,~y\rangle .
\end{eqnarray*}
Therefore $T^*x=PS^*x$ for every $x\in\h$. And so, for every $x\in\h$ we have
\begin{eqnarray*}
\|T^*x\|=\|PS^*x\|&\leq &\|S^*x\|=\|Sx\|=\|Tx\|~~~\mbox{because }\|P\|=1.
\end{eqnarray*}
$\hfill\blacksquare$\\
{\bf Question. }Does the converse of Theorem 3.1.7 holds?\\
The answer is negative. Let us consider the following operator $T=S^*+2S$ where $S$ is the unweighted shift in $\h$ i.e:
$Se_n=e_{n+1}$ for every $n\in\N$ with $(e_n)_n$ is an orthonormal basis of $\h$. A simple computation show that for every
$x=\sum\limits_n\alpha_ne_n\in\h$ we have
$$\big{\langle }(T^*T-TT^*)x~,~x\big{\rangle}=3\alpha_0^2\geq 0.$$
And so, from the equality (3.1) it follows that $\|T^*x\|\leq \|Tx\|$ for every $x\in\h$. On the other hand, for
$x=e_0-2e_2$ we have
$$\|(T^2)^{^*}x\|=\|T^{*^2}x\|=\sqrt{89}>\sqrt{80}=\|T^2x\|.$$
Hence by Theorem 3.1.7, $T^2$ is not subnormal operator. Therefore $T$ is not subnormal operator since every power of a
subnormal operator is subnormal.\\
So, it is natural to study this class of operators $T\in\lh$ which have the property
$$\|T^*x\|\leq \|Tx\|~~~\mbox{for every }x\in\h.$$
A such operator is called hyponormal.
\begin{prop}If $T$ is hyponormal operator on $\h$ then:\\
{\rm(i)}~~~$Tx=\lambda x$ implies $T^*x=\overline{\lambda}x$.\\
{\rm(ii)}~~If $Tx=\lambda x$ and $Ty=\mu y$ for $\lambda\not=\mu$ then $\langle x~,~y\rangle =0$.
\end{prop}
{\bf Proof. }The first property hold immediately since $T-\lambda I$ is also hyponormal operator. Now, for $\lambda\not=\mu$ we
have
\begin{eqnarray*}
\lambda\langle x~,~y\rangle &=&\langle\lambda x~,~y\rangle=\langle Tx~,~y\rangle\\
&=&\langle x~,~T^*y\rangle =\langle x~,~\overline\mu y\rangle~~~\mbox{by the last property}\\
&=&\mu \langle x~,~y\rangle .
\end{eqnarray*}
Since $\lambda\not=\mu$ then $\langle x~,~y\rangle =0$.$\hfill\blacksquare$
\begin{rem}
If $T\in\lh$ is a Hyponormal operator then for every $\lambda\in\C$, $\ker(T-\lambda I)$ is a closed invariant subspace for $T$
and $T^*$.
\end{rem}
\begin{thm}For every hyponormal operator $T\in\lh$,
$$\|T^n\|=\|T\|^n~~~~~~~~~\mbox{for every }n\geq 1.$$
\end{thm}
{\bf Proof. }For $n=1$, the equality is trivial; proceed by induction. For every vector $x\in\h$ we have
\begin{eqnarray*}
\|T^nx\|^2=\langle T^nx~,~T^nx\rangle&=&\langle T^*T^nx~,~T^{n-1}x\rangle\\
&\leq &\|T^*T^nx\|.\|T^{n-1}x\|\\
&\leq &\|T^{n+1}x\|.\|T^{n-1}x\|\\
&\leq &\|T^{n+1}\|.\|x\|.\|T^{n-1}\|.\|x\|=\|T^{n+1}\|.\|T^{n-1}\|.\|x\|^2.
\end{eqnarray*}
Since the vector $x$ is arbitrary, it follows that
$$\|T^n\|^2\leq \|T^{n+1}\|.\|T^{n-1}\|.$$
In view of the induction hypothesis ($\|T^k\|=\|T\|^k$ whenever $1\leq k\leq n$), this can be rewritten as
$$\|T\|^{2n}\leq \|T^{n+1}\|.\|T\|^{n-1}.$$
Therefore,
$$\|T\|^{n+1}\leq \|T^{n+1}\|.$$
Since the reverse inequality is universal, the induction step is accomplished.$\hfill\blacksquare$
\begin{cor}
If $T\in\lh$ is a subnormal operator and $S$ is its minimal normal extension on a Hilbert space ${\cal K}$ then $\|T\|=\|S\|.$
\end{cor}
{\bf Proof. }Since every subnormal operator is hyponormal (see Theorem 3.1.7), it follows from the last theorem that
$r(T)=\|T\|.$ On the other hand, we have $\sigma (S)\subset\sigma (T)$ (see Theorem 3.1.6) then
$$\|S\|=r(S)\leq r(T)=\|T\|.$$
The reverse inequality is trivial since $T$ is the restriction of $S$ on $\h$.$\hfill\blacksquare$
\begin{thm}Every hyponormal operator $T$ on $\h$ has the single valued extension property.
\end{thm}
{\bf Proof. }Let $F$ be a vector valued analytic function such that
$$(T-\lambda I)F(\lambda )=0~~\mbox{for every }\lambda\in D(F).$$
Then $$TF(\lambda )=\lambda F(\lambda )~~\mbox{for every }\lambda\in D(F).$$
Fix $\lambda\in D(F)$, then for every $\mu\not=\lambda\in D(F)$ we have by the last proposition
$$\langle F(\lambda )~,~F(\mu )\rangle =0.$$
By Pythagorean Theorem it follows that
$$\|F(\lambda )-F(\mu )\|^2=\| F(\lambda )\|^2+\|F(\mu )\|^2.$$
Letting $\mu\to\lambda$ gives $F(\lambda )=0.$
Therefore $F$ is identically zero since $\lambda$ is arbitrary element in $D(F)$. This complete the proof.$\hfill\blacksquare$
\begin{thm}[Stampfli and Radjabalipour]
Every hyponormal operator $T$ on $\h$ has the Dunford's Condition C (DCC).
\end{thm}
For the proof see \cite{21} and \cite{29}.$\hfill\blacksquare$
\section{Characterization of Subnormal Shifts}
\hspace*{1.00cm}The weighted shift operators are interesting for solving a lot of problems in operator theory, they can be used for examples and
counterexamples to illustrate many properties of operators. So, it is natural to ask which weighted shifts are normal, which are
hyponormal, and which are subnormal? The first question was already answered that there is no positive weight sequence that makes a
unilateral weighted shift normal and only the constant weight sequence that makes a bilateral weighted shift normal. The answer to
the second question is also easy that is the hyponormal weighted shifts are characterized by monotonically increasing weight sequences.
The answer to the third question is not easy. One formulation was offered by stampfli \cite{27}; another formulation due to C. Berger
is very different. It is elegant and easy to state \cite{16}, \cite{25}. \\
We begin by stating the {\bf Spectral Theorem for Normal Operators} (see \cite{10} and \cite{22}) which will be used throughout this section. A
vector $e\in\h$ is said to be star-cyclic for an operator $T\in\lh$ if $\h$ is the smallest closed invariant subspace for $T$ and
for its adjoint $T^*$ containing $e$. The operator $T$ is said to be star-cyclic if it has a star-cyclic vector. Note that a vector $e\in\h$ is
star-cyclic for a normal operator $T$ if and only if the closed linear vector subspace generated by
$\big{\{}T^nT^{*^m}e~:~n,~m\in\N\big{\}}$
coincide with $\h$. Every injective bilateral weighted shift is star-cyclic. For a positive measure $\mu$ with compact support
subset $K$ of \C, the linear map $N_{\mu}$ defined on $L^2-$space $L^2\big{(}K,\mu\big{)}$ by $N_{\mu}\big{(}f\big{)}(z)=zf(z)$ is
bounded operator on the Hilbert space $L^2\big{(}K,\mu\big{)}$; it is star-cyclic with star-cyclic vector the constant function $1$
since the set $C(K)$ of all complex continuous functions on $K$ is dense in $L^2\big{(}K,\mu\big{)}$.
\begin{thm}[Spectral Theorem for Normal Star-Cyclic Operators.] Let $e_0$ be a star-cyclic vector for a normal operator $T\in\lh$.
Then there is a probability measure $\mu$ in the closed disc $\overline {\D}$ of the complex plane $\C$ of radius $\|T\|$ such that
there is a unique isomorphism $V:\h\longrightarrow L^2\big{(}\overline {\D},\mu\big{)}$ with $Ve_0=1$ and $VTV^{-1}=N_{\mu}.$
\end{thm}
Let $T$ be a weighted shift on a Hilbert space $\h$ with a positive weight sequence $(\omega_n)_n$, that is
$$Te_n=\omega_n e_{n+1}=\frac{\beta_{n+1}}{\beta_n}e_{n+1}$$
where $(e_n)_n$ is a orthonormal basis of $\h$ and $\beta $ is the following sequence given
by:

\begin{displaymath}
\beta_n=\left\{
\begin{array}{ll}
\omega_0...\omega_{n-1}&\textrm{\mbox{if }$n>0$}\\
\\
1&\textrm{\mbox{if }$n=0$}\\
\\
\frac{1}{\omega_{n}...\omega_{-1}}&\textrm{\mbox{if }$n<0$}\\
\end{array}
\right.
\end{displaymath}

\begin{thm}[Berger]If $T$ is unilateral weighted shift then $T$ is subnormal if and only if there is a probability measure $\mu$
on the closed interval $[0,\|T\|]$ such that for $n\geq 1$
$$\beta_n^2=\int t^{2n}d\mu (t).$$
\end{thm}
{\bf Proof.} Without loss of generality, we assume that $\|T\|=1.$\\
Suppose first that $T$ is a subnormal operator and let $S$ be its minimal normal extension; then by the corollary 3.1.11 we have
$\|S\|=1$. Let $M$ be a closed invariant subspace for $S$ and its adjoint $S^*$ containing $e_0$, then the restriction
$S_{|M}$ of $S$ on $M$ is a normal operator; since $S^ne_0=T^ne_0=\beta_ne_n$ then $e_n\in M$ for every $n\geq 0$. Hence $\h\subset
M$; therefore by the minimality of $S$ it follows that $M={\cal K}$ and so, $e_0$ is star-cyclic vector for the normal operator $S$. The
Theorem 3.2.1 shows that there is a probability measure $\mu$ on the closed unit disc $\overline {\D}$ such that there is a unique
isomorphism $V:\h\longrightarrow L^2\big{(}\overline {\D},\mu\big{)}$ with $Ve_0=1$ and $VSV^{-1}=N_{\mu}.$ Let $\nu$ be the
probability measure in the closed unit interval $I$ defined for each Borel subset $E$ of $I$ by $\nu (E)=\mu \big{(}S^{-1}(E)\big{)}$
where $S:\overline {\D}\longrightarrow I$ is the mapping given by $S(z)=|z|.$ It follows that
$\int\limits_If(t)d\nu (t)=\int\limits_{\overline {\D}}f\circ S(z)d\mu (z)$ for every $L^1(I,\nu)-$function $f$. In particular, for
every $n\geq 0$,
we have
\begin{eqnarray*}
\int\limits_It^{2n}d\nu (t)&=&\int\limits_{\overline {\D}}|z|^{2n}d\mu (z)=\int\limits_{\overline {\D}}|N_{\mu}^n(1)|^2d\mu (z)\\
&=&\|N_{\mu}^n(1)\|^2=\|N_{\mu}^nVe_0\|^2=\|VS^ne_0\|^2\\   
&=&\|S^ne_0\|^2=\|T^ne_0\|^2\\
&=&\beta_n^2
\end{eqnarray*}
and so, the proof of the necessity of the Theorem is complete.\\
Now, let $\mu$ be a probability measure in the closed unit interval $I$ such that $\beta_n^2=\int t^{2n}d\mu (t)$ for every $n\geq
0.$
We consider the probability measure $\nu =d\mu (r)\frac{d\theta}{2\pi}$ on the closed unit disc $\overline {\D}$. The operator
$N_{\nu}$ of "multiplication by $z$" on the Hilbert space $L^2(\overline {\D},\nu )$ is a bounded normal operator. Let $H^2(\nu )$ be
the closed subspace spanned by $\{z^n~:~n\geq 0\}$. Then $H^2(\nu )$ is invariant closed subspace for $N_{\nu}$ and $N_{\nu}$
restricted to $H^2(\nu )$ is a subnormal operator $R\in{\cal L}\big{(}H^2(\nu )\big{)}$. On the other hand, for every $n,m\in\N$ we
have
\begin{eqnarray*}
\langle z^n~,~z^m\rangle &=&\int\limits_{\overline {\D}}r^{n+m}e^{i(n-m)\theta}d\nu(re^{i\theta})\\
&=&\frac{1}{2\pi}\int\limits_Ir^{n+m}dr\int\limits_0^{2\pi}e^{i(n-m)\theta}d\theta .
\end{eqnarray*}
Hence $(v_n)_{n\in\N}$ is orthonormal basis for $H^2(\nu )$ where $v_n=\frac{z^n}{\|z^n\|}=\frac{z^n}{\beta_n}$ for every $n\in\N$.
And so, for every $n\in\N$ we have
$$Rv_n=N_{\nu}\big{(}\frac{z^n}{\beta_n}\big{)}=\frac{1}{\beta_n}z^{n+1}=\omega_nv_{n+1}.$$
It follows that $T$ is a subnormal operator.$\hfill\blacksquare$
\begin{rem}
If $\mu$ is a probability measure on the interval $[0,\|T\|]$ such that for $n\geq 1$
$$\beta_n^2=\int t^{2n}d\mu (t),$$
then using the Cauchy-Schwartz Inequality we get
$$\beta_n^2\leq\beta_{n-1}\beta_{n+1}~~~~~~\mbox{for every }n\geq 1.$$
This condition is equivalent that the weight sequence $(\omega_n)_n$ is increasing that is the weighted shift $T$ is hyponormal.\\
\end{rem}
We close this section with the characterization of subnormal bilateral shifts.
\begin{thm}If $T$ is bilateral weighted shift then $T$ is subnormal if and only if there is a probability measure $\mu$
on the closed interval $[0,\|T\|]$ such that for $n\geq 1$ the functions $t^n,~t^{-n}\in L^1(\mu )$ and
\begin{eqnarray*}
~&&\beta_n^2=\int t^{2n}d\mu (t)\\
&&\frac{1}{\beta_{-n}^2}=\int t^{-2n}d\mu (t).
\end{eqnarray*}  
\end{thm}


\chapter{Bounded Point Evaluations for Cyclic Operators }
\section{Bounded Point Evaluations}
\hspace*{1.00cm}Throughout this section, $T$ will be a cyclic bounded operator on a Hilbert space $\h$ with cyclic vector $x$. A complex number
$\lambda\in\C$ is said to be a {\bf bounded point evaluation} of $T$ if there is a constant $M>0$ such that
$$|p(\lambda )|\leq M\|p(T)x\|$$
for every complex polynomial $p.$ The set of all bounded point evaluations of $T$ will be denoted by $B(T)$. Note that it follows
from the Riesz Representation Theorem (see \cite{6} and \cite{7}) that $\lambda\in B(T)$ if and only if there is a unique vector denoted $k_\lambda\in\h$ such that
$p(\lambda )=\langle p(T)x~,~k_\lambda\rangle$ for every complex polynomial $p$.
\begin{prop}$B(T)=\Gamma (T)$ the compression spectrum of $T$.
\end{prop}
{\bf Proof. }Let $\lambda\in B(T)$, then $x$ can not be in the closure of the range, $\overline{Im(T-\lambda I)}$, of $(T-\lambda I)$
otherwise $T^nx$ is in $\overline{Im(T-\lambda I)}$ for every $n\geq 0$ since $T$ and $(T-\lambda I)$ commutes; hence
$\overline{Im(T-\lambda I)}=\h $ and $k_\lambda=0$ which is impossible since $1=\langle x~,~k_\lambda\rangle$. 
Conversely,let $\lambda\in\Gamma (T)$ then $x$ is not in the closure of the range of $T-\lambda I$ since $x$ is also cyclic vector
for $T-\lambda I.$ Let $y\in\h$ be the orthogonal projection of $x$ onto the orthogonal complement of the range of $T-\lambda I$
and set $k_\lambda =\frac{1}{\langle x~,~y\rangle}y$. Since for every complex polynomial $p$ there is a complex polynomial $q$ such
that $p(z)=(z-\lambda )q(z)+p(\lambda ),$ then
\begin{eqnarray*}
\langle p(T)x~,~k_\lambda\rangle &=&\langle\big{[}(T-\lambda I)q(T)+p(\lambda )\big{]}x~,~k_\lambda\rangle\\
&=&\langle (T-\lambda I)q(T)x~,~k_\lambda\rangle +\langle p(\lambda )x~,~k_\lambda\rangle\\
&=&0+p(\lambda )\langle x~,~k_\lambda\rangle =p(\lambda ).
\end{eqnarray*}
$\hfill\blacksquare$  
\begin{prop}Let $\lambda\in\C$, the following statements are equivalent:\\
{\rm (i)}~~~$\lambda\in B(T).$\\
{\rm (ii)}~~$\ker \big{(}(T-\lambda)^*\big{)}$ is one dimensional.\\
{\rm (iii)}~$\overline\lambda\in\sigma_p(T^*).$
\end{prop}
{\bf Proof. }It is clear that (ii) implies (iii) and also (iii) and (i) are equivalent since
$\overline\Gamma(T)=\sigma_p(T^*)$ (see \cite{15}).\\
(i) implies (ii). It is clear that $k_\lambda\in\ker \big{(}(T-\lambda)^*\big{)}$. Let $u$ be a non zero
element of $\ker\big{(}(T-\lambda)^*\big{)}$, so it suffices to prove that there is $\alpha\not= 0\in\C$ such that $u=\alpha
k_\lambda$. 
Let $\alpha=\frac{1}{\overline{\langle u~,~k_\lambda\rangle}}$.  Since for every polynomial $p$ there is a polynomial
$q$ such that $p(z)=(z-\lambda )q(z)+p(\lambda ),$ then
\begin{eqnarray*}
\langle p(T)u~,~\alpha k_\lambda\rangle &=&\langle\big{[}(T-\lambda I)q(T)+p(\lambda )\big{]}u~,~\alpha k_\lambda\rangle\\
&=&\langle (T-\lambda I)q(T)u~,~\alpha k_\lambda\rangle +\langle p(\lambda )u~,~\alpha k_\lambda\rangle\\
&=&\langle q(T)u~,~(T-\lambda I)^*\alpha k_\lambda\rangle +p(\lambda )\langle u~,~\alpha k_\lambda\rangle\\
&=&p(\lambda).
\end{eqnarray*}
Hence, $u=\alpha k_\lambda$. Thus proves that $\ker \big{(}(T-\lambda)^*\big{)}$ is one dimensional. This finishes the
proof.$\hfill\blacksquare$\\ 
\begin{dft}A subset $O$ of $B(T)$ which is open in $\C$ is said to be an {\bf analytic set }for $T$ if for every $z\in\h$ the
function $\widehat z$ defined on $B(T)$ by $\widehat z(\lambda )=\langle z~,~k_\lambda\rangle$ is analytic on $O$. The largest
analytic set for $T$ is will denoted by $B_a(T)$ and every point of its will be called {\bf analytic bounded point evaluation} for $T$.  
\end{dft}
\begin{lem}A subset $O$ of $B(T)$ which is open in $\C$ is an analytic set for $T$ if and only if the function
$\lambda\longmapsto\|k_\lambda\|$ is bounded on compact subsets of $O$.
\end{lem}
{\bf Proof. }First suppose that $O$ is an analytic set for $T$ and $K$ is a compact subset of $O$. For every $z\in\h$ the
function $\widehat z$ is analytic on $O$, in particular, $\sup\limits_{\lambda\in K}|\langle z~,~k_\lambda\rangle|<+\infty$. So, it
follows from the Uniform Boundedness Principle that $$\sup\limits_{\lambda\in K}\|k_\lambda\|<+\infty .$$
Conversely, suppose that the function $\lambda\longmapsto\|k_\lambda\|$ is bounded on
compact subsets of $O$. Let $z\in\h$, then there is a sequence of polynomial $(p_n)_n$ such that
$\lim\limits_{n\to +\infty}\|p_n(T)x-z\|=0$. 
And so, for every compact subset $K$ of $O$ it follows by using Cauchy-Schwartz inequality that
$$\sup\limits_{\lambda\in K}|p_n(\lambda )-\widehat z(\lambda )|\leq\sup\limits_{\lambda\in K}\|k_\lambda\|.\|p_n(T)x-z\|.$$
Hence, $\widehat z$ is analytic function on $O$ (see \cite{}).$\hfill\blacksquare$\\
\begin{prop}[Williams] $\Gamma (T)\backslash\sigma_{ap}(T)\subset B_a(T)$.
\end{prop}
{\bf Proof. }Since $\sigma (T)=\Gamma (T)\cup\sigma_{ap}(T)$ then $\sigma (T)\backslash\sigma_{ap}(T)=\Gamma (T)\backslash\sigma_{ap}(T)$. On
the other hand, $\sigma (T)\backslash\sigma_{ap}(T)=int\big{(}\sigma (T)\big{)}\backslash\sigma_{ap}(T)$ since the  boundary of the spectrum
of $T$ is contained in $\sigma_{ap}(T)$. So, the set $O=\Gamma (T)\backslash\sigma_{ap}(T)$ is a subset of $B(T)$ which is open in $\C$. 
Let now $\lambda\in O$ then there is a positive constant $C$ such that $\|z\|\leq C\|(T-\lambda I)z\|$ for every $z\in\h$. So,
for $\mu\in\C$ such that $|\mu-\lambda |\leq \frac{1}{2C}$ we have 
\begin{eqnarray*}
\|z\|&\leq&C\|(T-\lambda I)z\|=C\|(T-\mu I)z+(\mu -\lambda)z\|\\
&\leq &C\|(T-\mu I)z\|+|\mu -\lambda|\|z\|\\
&\leq&C\|(T-\mu I)z\|+\frac{1}{2}\|z\|.
\end{eqnarray*}
Therefore $\|z\|\leq 2C\|(T-\mu I)z\|$ for each $z\in\h$ and for each complex number $\mu$ satisfying $|\mu-\lambda |\leq
\frac{1}{2C}$.
In particular, for every polynomial $p$ we have 
$$|p(\lambda )|\leq \|p(T)x\|\|k_\lambda\|\leq 2C\|(T-\mu I)p(T)x\|\|k_\lambda\|~~~\mbox{for } \mu\in\C,~|\mu-\lambda
|\leq\frac{1}{2C}.\leqno (*)$$
Now, let $\mu\in\C$ then for every polynomial $p$ there is a polynomial $q$ such that $p(t)=(t-\mu )q(t)+p(\mu ).$ Hence for
$\mu\in\C,~|\mu-\lambda |\leq\frac{1}{2C}$
\begin{eqnarray*}
|p(\mu )|&\leq &|p(\lambda )|+|\lambda -\mu ||q(\lambda )|\\
&\leq &\|p(T)x\|\|k_\lambda\|+2C|\lambda -\mu |\|(T-\mu I)q(T)x\|\|k_\lambda\|~~~~~~~~~~~~~~\mbox{by }(*)\\
&\leq &\|p(T)x\|\|k_\lambda\|+2C|\lambda -\mu |\|p(T)x-p(\mu )x\|\|k_\lambda\|\\
&\leq &\|p(T)x\|\|k_\lambda\|+2C|\lambda -\mu |\bigg{[}\|p(T)x\|+|p(\mu )|\|x\|\bigg{]}\|k_\lambda\|.
\end{eqnarray*}
So, if in addition, $\mu\in\C,~|\mu-\lambda |\leq M=\min\bigg{(}\frac{1}{2C}~,~\frac{1}{4C\|k_\lambda\|\|x\|}\bigg{)}$ then
\begin{eqnarray*}
|p(\mu )|&\leq &|\|p(T)x\|\|k_\lambda\|+\frac{1}{4C\|k_\lambda\|\|x\|}2C\bigg{[}\|p(T)x\|+|p(\mu )|\|x\|\bigg{]}\|k_\lambda\|\\
&\leq &\|p(T)x\|\|k_\lambda\|+\frac{1}{2\|x\|}\|p(T)x\|+\frac{1}{2}|p(\mu )|.
\end{eqnarray*}
It follows that $|p(\mu )|\leq \|p(T)x\|\bigg{(}2\|k_\lambda\|+\frac{1}{\|x\|}\bigg{)}$ for each $\mu\in\C,~|\mu-\lambda |\leq M$
and for every polynomial $p$. 
And so, $\|k_\mu\|\leq 2\|k_\lambda\|+\frac{1}{\|x\|}$ for every $\mu\in\C,~|\lambda -\mu |\leq M$ since the set
$\{p(T)x~:~p\mbox{ is a polynomial }\}$ is dense in $\h$ and every vector $k_\mu$ define a bounded linear functional
on $\h$ of norm $\|k_\mu\|$. Therefore, the function $\mu\longmapsto \|k_\mu\|$ is bounded on compact subsets of
$O$. Hence $O$ is an analytic set for $T$, therefore $\Gamma (T)\backslash\sigma_{ap}(T)\subset B_a(T).\hfill\blacksquare$\\
\\
\hspace*{1.00cm}Tavan Trent \cite{300} proved that the converse of this proposition holds for the operator $S_\mu$ of multiplication by $z$ on $H^2(\mu )$, the
closure of the polynomials in the $L^2(\mu)$ space where $\mu$ is a positive finite compactly supported Borel measure.
\begin{prop}[Tavan Trent]
Let $\mu$ be a positive finite compactly supported Borel measure. Then 
$$\Gamma (S_\mu)\backslash\sigma_{ap}(S_\mu)= B_a(S_\mu).$$
\end{prop}
{\bf Proof. }We first observe that the constant function $1$ is a cyclic vector for $S_\mu$ and $p(S_\mu )1=p$ for every polynomial
$p$. Since $B_a(S_\mu)\subset B(S_\mu)=\Gamma (S_\mu )$ (see Proposition 4.1.1) then it suffices to prove that
$B_a(S_\mu)\cap\sigma_{ap}(S_\mu )=\emptyset$. Now, let $\lambda\in B_a(S_\mu)$ then there is $r>0$ such that $B_a(S_\mu)$ contained a
closed disc $\D$ centered at $\lambda $ and of radius $r$ since $B_a(S_\mu)$ is open set. So, 
$$\sup\big{\{}\|k_\xi\|~:~\xi\in\D\big{\}}=C<\infty.$$  
Let $p$ be a polynomial then for $\gamma$ in the boundary of $\D$ i.e: $|\gamma -\lambda |=r$, we have $|(\gamma -\lambda )p(\gamma
)|\leq \|(S_\mu -\lambda )p\|\|k_\gamma\|$ since $\gamma$ is a bounded point evaluation for $S_\mu$. Hence,
$$|p(\gamma )|\leq \frac{C}{r}\|(S_\mu -\lambda )p\|~~~\mbox{for every }\gamma\in\C,~|\gamma -\lambda |=r.$$
Since the polynomial $p$ is analytic function on $\D$ then by the Maximum Modulus Principle,
$$|p(\gamma )|\leq \frac{C}{r}\|(S_\mu -\lambda )p\|~~~\mbox{for every }\gamma\in\D .$$
Hence,
\begin{eqnarray*}
\|p\|^2&=&\int\limits_{\D}|p|^2d\mu~+~\int\limits_{\C\backslash\D}|p|^2d\mu\\
&\leq &\bigg{(}\frac{C}{r}\bigg{)}^2\|(S_\mu -\lambda )p\|^2\mu (\D
)~+~\int\limits_{\C\backslash\D}|\frac{1}{z-\lambda}|.|{z-\lambda}p(z)|^2d\mu (z)\\
&\leq &\bigg{(}\frac{C}{r}\bigg{)}^2\|(S_\mu -\lambda )p\|^2\mu (\D 
)~+~\frac{1}{r^2}\int\limits_{\C\backslash\D}\bigg{|}\bigg{(}(S_\mu-\lambda)p\bigg{)}(z)\bigg{|}^2d\mu (z)\\
&\leq &\bigg{(}\frac{C}{r}\bigg{)}^2\|(S_\mu -\lambda )p\|^2\mu (\D
)~+~\frac{1}{r^2}\|(S_\mu-\lambda)p\|^2.
\end{eqnarray*}  
Therefore $\|p\|\leq M\|(S_\mu-\lambda)p\|$ for every polynomial $p$ where $M=\sqrt{\frac{C^2\mu (\D )~+~1}{r^2}}$. 
Since the polynomials are dense in $H^2(\mu )$ then $\|f\|\leq M\|(S_\mu-\lambda)f\|$ for every
$f\in H^2(\mu )$. That is $\lambda\not\in\sigma_{ap}(S_\mu )$. This complete the proof.$\hfill\blacksquare$\\
\hspace*{1.00cm}Several authors defined the concept of bounded point evaluation for a positive finite compactly supported Borel measure $\mu$ (see
\cite{9}, \cite{300},...)that is: $\lambda\in\C$ is a bounded point evaluation for $\mu$ if $\sup\bigg{[}\frac{|p(\lambda
)|}{\|p\|}\bigg{]}$
where the supremum is taken over all polynomials whose $L^2(\mu )$ norm is not zero. In fact, a point $\lambda\in\C$ is a bounded
point evaluation for $\mu$ if and only if it is a bounded point evaluation for the operator $S_\mu$ of multiplication by $z$ on
the space $H^2(\mu )$ (according to the above definition given by Williams \cite{32} for an arbitrary cyclic
operators). Williams followed exactly the same proof of Trent \cite{300} to prove Proposition 4.1.5 for an arbitrary cyclic
operators but he was neither able to exploit the analytic function theory used by Trent in the proof of the
Proposition 4.1.6 nor to prove the converse of Proposition 4.1.5. Thus he pose in \cite{32} the following question:\\ 
{\bf Question 1.} Is $B_a(T)=\Gamma (T)\backslash\sigma_{ap}(T)$?\\  
\\
\hspace*{1.00cm}Note that, an affirmative answer to Question 1 would show that $B_a(T)$ is independent of choice of the
cyclic vector for the operator $T$. Nevertheless, Williams proved this fact without knowing the answer to the
question 1. 
\begin{prop}$B_a(T)$ does not depend on the choice of the cyclic vector for $T$.
\end{prop}
{\bf Proof. }Let $y\in\h$ be also a cyclic vector for $T$. So for each $\lambda$ in $B(T)$ there is a unique vector $h_\lambda\in\h$
such that $p(\lambda)=\langle p(T)x~,~k_\lambda\rangle=\langle p(T)y~,~h_\lambda\rangle$ for every polynomial $p$.  
Therefore for every $\lambda$ in $B(T)$ there is $t_\lambda\not= 0\in\C$ such that $h_\lambda=t_\lambda k_\lambda$ since the
subspace
$\ker (T-\lambda I)^*$ is one dimensional (see Proposition 4.1.2). Now let $O$ be an analytic set for $T$ with respect to the cyclic vector
$x$, then the function $\lambda\longmapsto \|k_\lambda\|$ is bounded on compact subsets of $O$. We have 
$\widehat{y}(\lambda)=\langle y~,~k_\lambda\rangle=\langle y~,~\frac{1}{t_\lambda}k_\lambda\rangle=\frac{1}{\overline{t_\lambda}}$   
for each $\lambda\in B(T)$. It follows that the function $\lambda\longmapsto t_\lambda$ is continuous on $O$, and thus it is bounded on
compact subsets of $O$. Therefore the function $\lambda\longmapsto \|h_\lambda\|$ is bounded on compact subsets of $O$. It follows
from Lemma 4.1.4 that $O$ is an analytic set for $T$ with respect to $y$. By symmetry, $B_a(T)$ is independent of the choice of the
cyclic vector for $T$. $\hfill\blacksquare$\\
\\
\hspace*{1.00cm}Let $\h_1$ and $\h_2$ be two Hilbert spaces. Two operators $R\in{\cal L}(\h_1)$ and $S\in{\cal L}(\h_2)$ are said to be {\bf
quasisimilar} if there exist two bounded transformations $X:\h_1\longrightarrow\h_2$ and $Y:\h_2\longrightarrow\h_1$
having trivial kernels and dense ranges such that $XR=SX$ and $RY=YS$. M. Raphael \cite{23} showed that quasisimilar cyclic subnormal 
operators have the same bounded point evaluations and the same analytic bounded point evaluations. The following theorem \cite{32} shows  
in general that two quasisimilar cyclic operators have the same bounded point evaluations and the same analytic bounded point evaluations. 
\begin{thm}Let $\h_1$ and $\h_2$ be two Hilbert spaces. Two quasisimilar cyclic operators $T\in{\cal L}(\h _1)$ and $S\in{\cal
L}(\h _2)$ have the same bounded point evaluations and the same analytic bounded point evaluations, i.e: $B(T)=B(S)$ and $B_a(T)=B_a(S)$.
\end{thm} 
{\bf Proof. }Suppose that there exist two bounded transformations $X:\h _1\longrightarrow\h _2$ and $Y:\h _2\longrightarrow\h _1$ having trivial kernels 
and dense ranges such that $XR=SX$ and $RY=YS$.
Since $X$ has dense range then $0\not\in\Gamma(X)$ so, $X^*$ is injective since $\overline\Gamma(X)=\sigma_p(X^*)$. If
$\lambda\in\C$ is eigenvalue of $S^*$ with the corresponding eigenvector $x\in\h_2$ then $R^*(X^*x)=\lambda X^*x$, and so,
$\lambda\in\sigma_p(R^*)$ since $X^*x\not=0$. Hence, $\sigma_p(S^*)\subset\sigma_p(R^*)$. By symmetry, $\sigma_p(S^*)=\sigma_p(R^*)$, thus
$B(R)=B(S)$. Suppose that $v\in\h_2$ is a cyclic vector for $S$ then $u=Yv$ is a cyclic vector for $R$. Let $\lambda$ be in $ B(R)=B(S)$
there are vectors $k_\lambda\in \h_1$ and $h_\lambda\in \h_2$ such that $p(\lambda)=\langle p(R)u~,~k_\lambda\rangle=\langle
p(S)v~,~h_\lambda\rangle$ for each polynomial. On the other hand $p(R)Y=Yp(S)$ since $RY=YS$ for every polynomial $p$. Ans so,
\begin{eqnarray*}
\langle p(S)v~,~Y^*k_\lambda\rangle &=&\langle Yp(S)v~,~k_\lambda\rangle\\
 &=&\langle p(R)Yv~,~k_\lambda\rangle\\
&=&\langle p(R)u~,~k_\lambda\rangle\\
&=&p(\lambda)\\
&=&\langle p(S)v~,~h_\lambda\rangle
\end{eqnarray*}
for every polynomial $p$ and $\lambda\in B(R)=B(S)$.
Hence, $h_\lambda=Y^*k_\lambda$ for every $\lambda\in B(R)=B(S)$. Since $\lambda\longmapsto \|k_\lambda\|$ is bounded on compact subsets of
$B_a(R)$ then $\lambda\longmapsto \|h_\lambda\|$ is bounded on compact subsets of $B_a(R)$. It follows from Lemma 4.1.4 that $B_a(R)\subset
B_a(S).$ By symmetry, the desired result holds i.e: $B_a(R)=B_a(S).\hfill\blacksquare$
\section{Analytic Bounded Point Evaluations for Unilateral Weighted Shift}
\hspace*{1.00cm}In \cite{25}, A.L. Shields represented a weighted shift operator as ordinary shift operator (that is, as "multiplication by z")
on a Hilbert space of formal power series (in the unilateral case) or formal Laurent series (in the bilateral case). Thus
he defined the concept of bounded point evaluations of a weighted shift to examine which power series and Laurent series
represent analytic functions. In fact, this concept of bounded point evaluations for injective unilateral weighted shift 
coincide with the one defined by L.R. Williams in \cite{32}.\\  
\\
\hspace*{1.00cm}We now describe the set of bounded point evaluations and the set of analytic bounded point evaluations for an arbitrary injective unilateral
weighted shift. Let $S$ be a unilateral weighted shift on a Hilbert space $\h$ with a positive weight sequence $(\omega_n)_{n\geq 0}$, that is
$$Se_n=\omega_n e_{n+1}=\frac{\beta_{n+1}}{\beta_n}e_{n+1}$$
where $(e_n)_{n\geq 0}$ is a orthonormal basis of $\h$ and $\beta $ is the following sequence given
by:
\begin{displaymath}
\beta_n=\left\{
\begin{array}{ll}
\omega_0...\omega_{n-1}&\textrm{\mbox{if }$n>0$}\\
\\
1&\textrm{\mbox{if }$n=0$}
\end{array}
\right.
\end{displaymath}

The unilateral weighted shift $S$ is cyclic of cyclic vector $e_0$. It follows from Corollary 1.3.6 and Theorem 4.1.8 that the set of the bounded point
evaluations and the set of analytic bounded point evaluations for the unilateral weighted shift $S$ have a circular symmetry about the origin.
Also, it follows from Theorem 1.3.12 and Proposition 4.1.1 that $B(S)=\{0\}$ if $r_2(S)=0$ otherwise 
$$\big{\{}\lambda\in\C~:~|\lambda |<r_2(S)\big{\}}\subset B(S)\subset\big{\{}\lambda\in\C~:~|\lambda
|\leq r_2(S)\big{\}},$$
where $r_2(S)=\liminf\limits_{{n\to +\infty}}(\beta_n)^{\frac{1}{n}}.$
\begin{thm}If $r_2(S)>0$ then $B_a(S)=\big{\{}\lambda\in\C~:~|\lambda |<r_2(S)\big{\}}$. 
\end{thm} 
{\bf Proof. }Let $\lambda\in B(S)$ then there is $k_\lambda=\sum\limits_{n\geq 0}\widehat{k_\lambda}(n)e_n\in\h$ such that $p(\lambda)=\langle
p(S)e_0~,~k_\lambda\rangle$ for every polynomial $p$. In particular, for every $n\geq 0$ we have
\begin{eqnarray*}
\widehat{k_\lambda}(n)&=&\langle e_n~,~k_\lambda\rangle\\
&=&\langle \frac{1}{\beta_n}S^ne_0~,~k_\lambda\rangle\\
&=&\frac{\lambda^n}{\beta_n}.
\end{eqnarray*} 
Now, let $h=\sum\limits_{n\geq 0}a_ne_n$ then $\widehat h(\lambda)=\langle h~,~k_\lambda\rangle=\sum\limits_{n\geq
0}\frac{a_n}{\beta_n}\lambda^n$. In fact this series is absolutely convergent since $\sum\limits_{n\geq 0}|a_n|e_n$ is also in $\h$ and
$|\lambda|$ is also in $B(S)$. Hence, $\widehat h$ is analytic in the interior of $B(S)$ which is exactly the disc
$\big{\{}\lambda\in\C~:~|\lambda |<r_2(S)\big{\}}$. This proves the theorem. $\hfill\blacksquare$\\
\hspace*{1.00cm}For the weighted shift $S$, the spectrum $\sigma(S)$ is known to be the disk $\big{\{}\lambda\in\C~:~|\lambda|\leq r(S)\big{\}}$, 
(see Theorem 1.3.8) and the approximate point spectrum $\sigma_{ap}(S)$ is known to be the annulus $\big{\{}\lambda\in\C~:~r_1(S)\leq|\lambda|\leq r(S)\big{\}}$
(see Theorem 1.3.10).
Therefore, $\Gamma(S)\backslash\sigma_{ap}(S)=\sigma(S)\backslash\sigma_{ap}(S)=\big{\{}\lambda\in\C~:~|\lambda|< r_1(S)\big{\}}$.
And so, $\Gamma(S)\backslash\sigma_{ap}(S)\subsetneqq B_a(S)$ if and only if $r_1(S)<r_2(S)$. Hence, a negative answer
to the question 1 (see \cite{32}) can be given by a unilateral weighted shifts $S$ for which $r_1(S)<r_2(S)$. Let us consider
an example of a such weighted shift.
For $s\in\N$ there are unique $n,~k\in\N$ such that $s=n!+k$ with $0\leq k\leq (n+1)!-n!-1$. We set $$\beta_s=\beta_{n!+k}=e^k.$$
 And so, for every $k\in\N$, we have
\begin{displaymath}
\frac{\beta_{k+1}}{\beta_k}=\left\{
\begin{array}{ll}
e &\textrm{\mbox{if }$n!\leq k<(n+1)!-1$}\\
\\
\frac{1}{e^{(n+1)!-n!-1}}&\textrm{\mbox{if }$k=(n+1)!-1$}
\end{array}
\right.
\end{displaymath}
Hence, the unilateral weighted shift $S$ corresponding to the weight $\big{(}\frac{\beta_{n+1}}{\beta_n}\big{)}_{n\geq 0}$ is bounded (see Proposition 1.3.1).
For every $n\geq 2$, set $k_n=n!-n$. Clearly, we have $(n-1)!\leq k_n<n!$ and $\beta_{k_n}=e^{\big{(}n!-(n-1)!-n\big{)}}$.
And so,
$$\inf\limits_{k\geq 0}\frac{\beta_{n+k}}{\beta_k}\leq\frac{\beta_{n+k_n}}{\beta_{k_n}}=
\frac{1}{e^{\big{(}n!-(n-1)!-n\big{)}}}.$$
Hence,  
$$r_1(S)=\lim\limits_{n\to\infty}\bigg{[}\inf\limits_{k\geq 0}\frac{\beta_{n+k}}{\beta_k}\bigg{]}^{\frac{1}{n}}=0.$$
On the other hand, it is clear that
$$r_2(S)=\liminf\limits_{n\to\infty}\bigg{[}\beta_n\bigg{]}^{\frac{1}{n}}=1.$$
Therefore, $B_a(S)=\big{\{}\lambda\in\C~:~|\lambda|<1\big{\}}$ and $\Gamma(S)\backslash\sigma_{ap}(S)=\emptyset.$\\  
For the unilateral weighted shift $S$ the set of its analytic bounded point evaluations is exactly the interior of the set of its bounded
point evaluations. This suggests the following question, which was first posed by J.B.Conway in \cite{9} page 65:\\
{\bf Question 2:} Is always the interior of $B(T)$ coincide with $B_a(T)$ for an arbitrary cyclic operator $T\in\lh$?\\

\section{The Local Spectra Through Bounded Point Evaluations}
\hspace*{1.00cm}In what follows, $\h$ denote a Hilbert space. In this section we will show that for a cyclic operator $T\in\lh$
with Dunford's Condition C and without point spectrum, the local spectra of $T$ at $x$ is equal to the spectrum of $T$ for each $x$ in a dense subset of $\h$. 
\begin{lem}Suppose that $\h=\h_1\oplus\h_2$ where $\h_1$ and $\h_2$ are two Hilbert spaces such that $\h_2$ is finite dimensional. If
$T\in\lh$ has the single valued extension property and $\h_1$ is invariant subspace for $T$ then $A=T_{|\h_1}$ has the single valued
extension property and $\sigma_{_A}(x)=\sigma_{_T}(x)$ for every $x\in\h_1$.
\end{lem}
{\bf Proof. }It is clear that $A$  has the single valued extension property and $\sigma_{_T}(x)\subset\sigma_{_A}(x)$ for every
$x\in\h_1$. Conversely, let $x\in\h_1$ then $\widetilde x=F_1+F_2$ on $\rho_{_T}(x)$ where 
$F_i=P_i\widetilde x$ and $P_i:\h\longrightarrow\h_i$ are the orthogonal projections $(i=1,~2)$. And so, for every $\lambda\in\rho_{_T}(x)$
\begin{eqnarray*}
x&=&(T-\lambda I)\widetilde x(\lambda)\\
&=&(T-\lambda I)F_1(\lambda)+(T-\lambda I)F_2(\lambda)\\
&=&(A-\lambda I)F_1(\lambda)+P_1(T-\lambda I)F_2(\lambda)+P_2(T-\lambda I)F_2(\lambda)\\
&=&(A-\lambda I)F_1(\lambda)+P_1TF_2(\lambda)+(P_2T-\lambda I)F_2(\lambda)\\
&=&\bigg{[}(A-\lambda I)F_1(\lambda)+P_1TF_2(\lambda)\bigg{]}+(P_2T-\lambda I)F_2(\lambda).
\end{eqnarray*}  
It follows that $(P_2T-\lambda I)F_2(\lambda)=0$ for every $\lambda\in\rho_{_T}(x)$. Since $\sigma(P_2T)$ is a finite set then $F_2$ is
identically zero function. Hence,
$$x=(T-\lambda I)\widetilde x(\lambda)=(A-\lambda I)F_1(\lambda)~~~~\mbox{for all }\lambda\in\rho_{_T}(x).$$
Therefore, $\sigma_{_A}(x)\subset\sigma_{_T}(x)$. The proof is complete.$\hfill\blacksquare$\\
\begin{lem}
Suppose that $\h=\h_1\oplus\h_2$ where $\h_1$ and $\h_2$ are two Hilbert spaces such that $\h_2$ is finite dimensional. If $\h_1$ is an
invariant subspace for an operator $T\in\lh$ which satisfies DCC then the operator $A=T_{|\h_1}$ satisfies the DCC.
\end{lem}
{\bf Proof. }It follows from the last Lemma that for every closed subset $F$ of $\C$, ${\cal H}_{1_{_A}}(F)={\cal H}_{_T}(F)\cap\h_1$. So, the desired
result holds.$\hfill\blacksquare$\\
\begin{lem}
Suppose that $\h=\h_1\oplus\h_2$ where $\h_1$ and $\h_2$ are two Hilbert spaces such that $\h_2$ is finite dimensional. If
$T\in\lh$ is a operator without point spectrum such that $\h_1$ is invariant subspace for $T$ then $\sigma(T)=\sigma(A)$ where $A=T_{|\h_1}$.
\end{lem}
{\bf Proof. }We first observe that there are two operators
$B:\h_2\longrightarrow\h_1$ and $C:\h_2\longrightarrow\h_2$ such that for every $x=x_1\oplus x_2\in\h$,
$Tx=\big{[}Ax_1+Bx_2\big{]}+Cx_2$. Now, suppose that $A$ is invertible in ${\cal L}(\h_1)$ then $T\h_1=\h_1$ and $\h_1\cap T\h_2=\{0\}$.    
Let $x_2\in\h_2$ such that $Cx_2=0$ then $Tx_2=Bx_2\in\h_1\cap T\h_2$, hence $C$ is injective
and so, by the finite-dimensionality, $C$ is invertible. Therefore
the operator $T$ is invertible with inverse given by: 
$$T^{-1}x=\big{[}A^{-1}x_1-A^{-1}BC^{-1}x_2\big{]}+C^{-1}x_2~~~\mbox{for every }x=x_1\oplus x_2\in\h.$$
Thus $\sigma(T)\subset\sigma(A)$. The converse follows from Theorem 2.1.2 and Lemma 4.3.1. $\hfill\blacksquare$\\  
\begin{thm}
Let $T\in\lh$ be a cyclic operator with cyclic vector $x\in\h$ and let $S\in\lh$ be an operator commutes with $T$ such that
$\ker (S^*)$ is finite dimensional. If $T$ satisfies DCC and $\sigma_p(T)=\emptyset$ then $\sigma_{_T}(Sx)=\sigma(T).$  
\end{thm}  
{\bf Proof. }Let $\h_1$ be the closed linear subspace generated by $\big{\{}T^nSx~:~n\geq 0\big{\}}$ then it is clear that $\h_1$ is invariant
subspace for $T$. Since $TS=ST$ then for every $n\geq 0$, $T^nx\in\h_1$; therefore $\h_1=\overline{Im(S)}$ since $x$ is a cyclic
vector for $T$. Thus $\h=\h_1\oplus\h_2$ where $\h_2=\ker(S^*)$. It follows from the Lemma 4.3.1 that $\sigma_{_T}(Sx)=\sigma_{_A}(Sx)$ where 
$A=T_{|\h_1}$.  
Since $A$ is a cyclic operator with cyclic vector $Sx$ and satisfies DCC (Lemma 4.3.2) then it follows from Proposition 2.2.2 that
$\sigma_{_A}(Sx)=\sigma(A).$  
Since $\sigma_p(T)=\emptyset$, it follow from Lemma 4.3.3 that $\sigma(T)=\sigma(A)$. And so, 
$\sigma(T)=\sigma(A)=\sigma_{_A}(Sx)=\sigma_{_T}(Sx).$ 
The proof is complete.$\hfill\blacksquare$
\begin{rem}
Let $T\in\lh$ be a cyclic operator with cyclic vector $x\in\h$ and satisfies DCC such that $\sigma_p(T)=\emptyset$. 
Since for every non zero polynomial $p$, there exists $a,\alpha_1,\ldots,\alpha_n\in\C$ such that
$p(T)^*=a(T^*-\alpha_1 I)\ldots(T^*-\alpha_n)$. Then it follows from Proposition 4.1.2 that $\ker p(T)^*$ is finite
dimensional. Hence, $\sigma_{_T}\big{(}p(T)x\big{)}=\sigma(T)$ for every
non zero polynomial $p ~(see ~{\rm Theorem ~4.3.4})$. Therefore, $\sigma_{_T}(y)=\sigma(T)$ holds for all non zero $y$ in a dense subset of $\h$. 
{\rm L.R. Williams} proved in {\rm Theorem 2.5} of {\rm\cite{33}} that if $T$ is a non normal hyponormal (unilateral or bilateral) 
weighted shift operator then $\sigma_{_T}(x)=\sigma(T)$ for every non zero element $x\in\h$. 
\end{rem}

We state and we give a simple proof of Theorem 2.5 of \cite{33} using the fact that a non zero analytic function has isolate zeroes.
\begin{thm}Let $T$ be a non normal hyponormal weighted shift on $\h$. Then for every a non zero element $x\in\h$, 
$\sigma_{_T}(x)=\sigma(T)$.
\end{thm}
{\bf Proof. }First suppose that $T$ is a non normal hyponormal unilateral weighted shift, then $r(T)=r_1(T)=r_2(T)=\|T\|>0$. 
Let now $x\in\h$ such that there is $\lambda\in\sigma(T)\backslash\sigma_{_T}(x)$. So, there is vector valued analytic
function $F$ on an open neighbourhood $V$ of $\lambda$ such that
$$(T-\mu I)F(\mu)=x~~\mbox{for every }\mu\in V.$$
Since $\emptyset\not= V\cap int\big{(}\sigma(T)\big{)}\subset B_a(T)$, then for every $\mu\in V\cap int\big{(}\sigma(T)\big{)}$ we have,
\begin{eqnarray*}
\widehat{x}(\mu)&=&\langle x~,~k_\mu\rangle\\
&=&\langle (T-\mu I)F(\mu)~,~k_\mu\rangle\\
&=&\langle F(\mu)~,~(T-\mu I)^*k_\mu\rangle\\
&=&0.
\end{eqnarray*}
Hence, $\widehat{x}\equiv 0$. And so, $x=0$.\\
\hspace*{1.00cm}The case of a non normal hyponormal bilateral weighted shift is similar using the bounded point evaluations of the sense
of A.L. Sheilds \cite{25}.
\begin{rem}Using the same proof of this Theorem, one can see that for every injective unilateral weighted shift $T$, $\overline{B_a(T)}\subset \sigma_{_T}(x)$ for every 
non zero element $x\in\h$. The same for every bilateral weighted shift.
\end{rem}
\noindent{\bf Question 3.} For which hyponormal operators $T$ we have $\sigma_{_T}(h)=\sigma(T)~\forall h\not= 0$?
\chapter*{Acknowledgements}
\addcontentsline{toc}{chapter}{Acknowledgements}
I humbly acknowledge the blessings of {\bf Almighty Allah} who gave me this opportunity to join 
the ICTP Diploma Program Mathematics and who has enabled me to complete my dissertation.\\
\\
I would like to express my deepest gratitude to my supervisor {\bf Professor C.E. Chidume}, 
Co-ordinator Diploma Course Mathematics. His supervision, his help and his moral support made learning from him 
very pleasant and valuable experience. Next, I am extremely grateful to {\bf Professor E.H. Zerouali}, Regular 
Associate at ICTP, for his help and continuous guidance in finishing my dissertation.\\
\\
Many thanks goes to {\bf Ms. Concetta Mosca} Diploma Office Incharge for her kind and countinous help and useful advices throughout the Diploma 
Program. \\
\\
Finally, I am thankful to {\bf UNESCO}, the {\bf IAEA} and {\bf Professor M.A. Virasoro} (Director of ICTP) for their kind hospitality
at the centre during the Diploma Program.\\
\\  
\\
\\
Trieste Italy, August 2000 \hfill{\bf A. BOURHIM}


\end{document}